\documentclass[3p,super,authoryear]{elsarticle}

\usepackage{amssymb,amsmath,mathtools}
\usepackage{amsthm}
\usepackage{xcolor}
\usepackage{stmaryrd}
\usepackage{multirow}
\usepackage{graphbox}
\usepackage[hyphens]{url}
\usepackage[colorlinks=true]{hyperref}
\renewcommand{\vec}[1]{\boldsymbol#1}
\renewcommand{\Vec}[1]{\mathbf#1}
\newcommand{\tensorTwo}[1]{\boldsymbol#1}
\newcommand{\tensorFour}[1]{\mathbb{#1}}
\newcommand{\af}[1]{\textcolor{black}{#1}}
\newcommand{\pt}[1]{\textcolor{black}{#1}}
\newcommand{\oneway}[1]{\textcolor{black}{#1}}
\usepackage[scaled]{helvet}
 
\usepackage[T1]{fontenc}

\usepackage{titlesec}
\usepackage{tikz}
\usepackage{xcolor}

\newcommand{\whiteline}[1]{%
  \begin{tikzpicture}[baseline=(text.base)]
    \node[inner sep=0pt, outer sep=0pt] (text) {\strut #1};
    \draw[white, line width=2pt] (text.south west) -- (text.south east);
  \end{tikzpicture}%
}

\titleformat{\section}
  {\normalfont\large\bfseries}
  {}
  {0pt}
  {\whiteline}

\titleformat{\subsection}
  {\normalfont\bfseries\itshape}
  {}
  {0pt}
  {\whiteline}

\usepackage[nomarkers,noheads]{endfloat}

\journal{Netherlands Journal of Geosciences}

\begin{document}

\begin{frontmatter}



\title{Unexpected fault activation due to underground gas storage in produced reservoirs. Part I: 
Mathematical model and mechanisms}


%
%
%
%
%
%
%

\author[UNIPD]{Andrea Franceschini\corref{correspondingauthor}}
\ead{andrea.franceschini@unipd.it}

\author[UNIPD]{Selena Baldan}
\ead{selena.baldan@phd.unipd.it}

\author[UNIPD]{Massimiliano Ferronato}
\ead{massimiliano.ferronato@unipd.it}

\author[UNIPD]{Carlo Janna}
\ead{carlo.janna@unipd.it}

\author[UNIPD]{Claudia Zoccarato}
\ead{claudia.zoccarato@unipd.it}

\author[M3E]{Matteo Frigo}
\ead{m.frigo@m3eweb.it}

\author[M3E]{Giovanni Isotton}
\ead{g.isotton@m3eweb.it}

\author[UNIPD]{Pietro Teatini}
\ead{pietro.teatini@unipd.it}

\affiliation[UNIPD]{organization={Department of Civil, Environmental and Architectural
                                  Engineering, University of Padova},
                    city={Padova},
                    postalcode={35131},
                    country={Italy}}
\affiliation[M3E]{organization={M3E S.r.l.},
                  city={Padova},
                  postalcode={35121},
                  country={Italy}}
\cortext[correspondingauthor]{Corresponding author}

\begin{abstract}

\pt{
Underground gas storage (UGS) is a critical technology for managing seasonal gas consumption peaks, increasingly important in the face of market uncertainties. However, safety concerns arise when reactivating pre-existing faults in faulted basins, where human activities may trigger seismic events. Typically, faults are reactivated when shear stress exceeds a critical frictional threshold, but unexpected fault reactivations have been observed during cushion gas injection (CGI) and UGS cycles in the Netherlands, even when the stress regime suggests stability. This two-part study introduces a novel simulation framework to better understand the mechanisms behind fault reactivation in complex settings such as the Rotliegend formation in the Netherlands. 
A 3D mathematical model coupling frictional contact mechanics in faulted porous rocks with fluid flow allows for predictive analysis of fault behavior. The effect of the storage of different fluids for various purposes, such as the long-term sequestration of CO$_2$, the regular injection and extraction cycles of CH$_4$, and the highly irregular cycles of H$_2$, is investigated with respect to fault activation hazard. The ultimate goal is to define a safe operational bandwidth for UGS activities in faulted reservoirs.}

\pt{
Part I of this work presents this comprehensive simulation tool where a slip-weakening constitutive law is introduced to model fault behavior. The approach is designed to address the complex geological setting that characterizes the Rotliegend formation, where multiple factors influence the behavior of fault systems. We succeed in explaining and modeling the occurrence of unexpected fault reactivations. The analysis shows that reactivation during primary production (PP) causes stress redistribution, leading to a new deformed equilibrated configuration. When the fault is reloaded in the opposite direction during cushion gas injection (CGI) or UGS cycles, further activation can occur, even if the stress does not exceed the initial or maximum stress previously experienced by the formation.}

\end{abstract}




\end{frontmatter}


\section{Introduction}
Seismicity associated with fluid withdrawal from and injection into deep reservoirs is a
geomechanical hazard that is receiving increasing attention in the scientific literature
\cite{Dog18, Fou_etal18, KerWei18}. Reactivation of faults, both aseismic and seismic, is
caused by the change in the natural stress regime on the discontinuity surface due to variations in the
pore pressure $p$ in the reservoir where mining activities are carried out. More
specifically, the onset and amount of slip, and the size of the reactivated fault zone
depend on how the stress changes caused by human operations at depth can interfere
with the natural stress regime \cite{Seg_etal94, Het_etal00, Can_etal19}.

\begin{figure}
  \centering
  \includegraphics[width=0.75\textwidth]{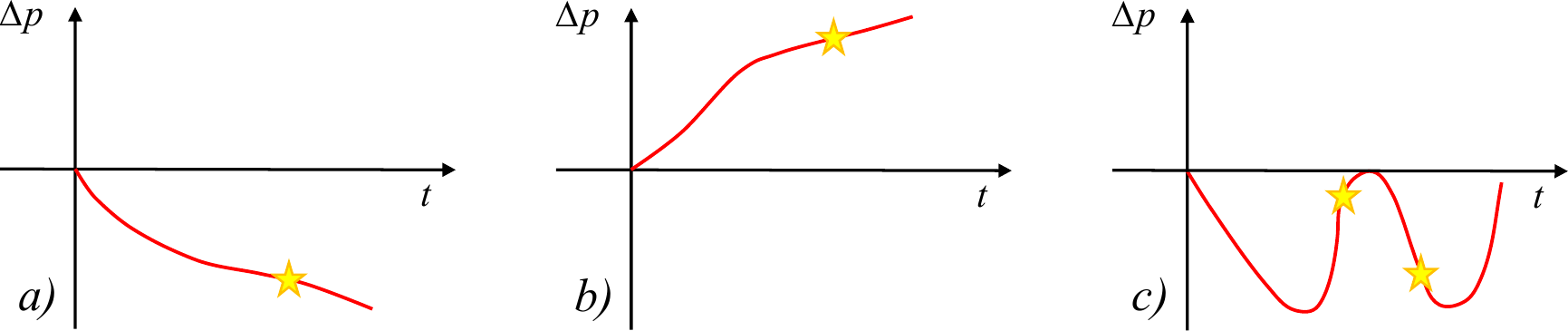}
  \caption{Sketches of two (\textit{a} and \textit{b}) ``expected'' induced seismicity
scenarios and one (\textit{c}) ``unexpected''. \textit{a}) Primary production with large
pressure drop, \textit{b}) fluid injection (CO$_2$ sequestration, waste water disposal,
fracking) with significant pressure increase), \textit{c}) pressure in the range already
experienced (UGS with $p < p_i$).}
  \label{fig:sketch}
\end{figure}

Current state-of-the-art research on this topic focuses on two main processes: i)
seismicity induced by production of (conventional) hydrocarbon reservoirs, where pore
pressure depletion $\Delta p$ and differential reservoir compaction are the main factors
generating fault reactivation \cite{Seg_etal94, Bui_etal17} (Fig.~\ref{fig:sketch}a); ii)
fluid injection at depth (CO$_2$ sequestration, production from unconventional reservoirs,
enhanced geothermal systems) where,
even without considering possible thermal processes,
the increase in fluid pressure (largely) above the natural undisturbed value $p_i$ within
fault zones that cross or bound the targeted formation drives the reactivation of
rock discontinuities (Fig.~\ref{fig:sketch}b) \cite{ChaSeg16, Rut_etal16, Eyr_etal19,
Yuy_etal19}.

In the last decade, induced seismicity has been observed in some parts of the world, also
in reservoirs used for underground gas storage (UGS). Somehow unexpectedly, fault
reactivation occurred not only during primary production (PP) or gas storage at pressure
larger than $p_i$ \cite{Def_etal95, Ces_etal14, Zho_etal19}, i.e., at a stress regime that
had never been experienced before by the reservoir and the nearby faults, but also during
cushion gas injection (GCI) or production and storage phases with a pore pressure smaller
than $p_i$ and larger than $p_{\mathrm{min}}$, i.e., the minimum pressure experienced by
the field usually at the end of primary production before its conversion to UGS
(Fig.~\ref{fig:sketch}c) \cite{MIT09, Kra_etal13, TNO15, NAM16}.

The purpose of this work is to shed light on these ``unexpected'' events. Because of the
current CH$_4$ importance for energy production purposes and the international turbulence
on this market, interest in developing UGS projects is increasing worldwide. Currently, multiple elements characterize UGS: seasonal and short-term balancing, strategic reserves
in case of interruption of deliveries, optimization of gas production and gas system
distribution, overcoming local restrictions of gas grids \cite{Pla09, Cor19}. More
recently, UGS has also been investigated as a possible method of storing green energy in
terms of compressed air and H$_2$ \cite{ZakSyr15, Mou_etal19, Sop_etal19}. Sources of
green energy, such as wind, waves, and sun, are characterized by a natural high-frequency
fluctuation (from hours to day/night and to weeks). Excess electricity can be used to
synthesize hydrogen or compress air, store the gas in deep aquifers or depleted
reservoirs, and use it at a later stage as fuel to generate electricity.
The same technology can also be applied for long-term geological sequestration activities,
such as CO$_2$ capture and sequestration to reduce carbon dioxide emissions in the
atmosphere \cite{MosMedBabTur19}. In this case, the targeted reservoir does not undergo a
cyclic loading/unloading strength, but pressure can increase to a steady state
value, usually smaller or equal to the initial value $p_i$. 

Analysis of the social and
environmental hazards associated with subsurface gas storage is a recurring
issue whenever a new UGS site is planned. Many different aspects are involved, such as
formation integrity, health and safety in relation to public perception, economic hazard, and
environmental impact. Among the latter, the geomechanical effects induced by seasonal gas
injection and withdrawal, such as movements of the land surface, may play an important
role \cite{Tea_etal11}.

UGS has rarely been associated with induced seismicity. According to data provided by the
European Commission \cite{EurComm}, recent works \cite{li2022national}, and the
\textit{HiQuake} database \cite{HiQuake}, only a few sites have reported human-induced
earthquakes out of 160 UGS facilities in Europe and more than 380 in the USA
\cite{Fou_etal18}. Three cases of these, i.e., Bergermeer, Norg, and Grijpskerk fields,
are located in The Netherlands.
These reservoirs are located in the Carboniferous-Rotliegend formation, northern Europe,
which is one of the most intensively explored petroleum systems in the world \cite{Gau03},
and where a relatively high number of induced seismicity events have been recorded over the
last decades \cite{Wee_etal14, Uta17}. Several studies have addressed the topic of fault
reactivation in Rotliegend reservoirs, the most famous of which concerns the Groningen
field \cite{ThiBre15}. Most studies focus on a specific reservoir in the Netherlands and
northern Germany, or, more generally, try to investigate the relationship between the
typical geological features of these reservoirs, their usual production life and the
possible induced seismicity \cite{Bui_etal17, Wee_etal17, Zbi_etal17, Hau_etal18,
Can_etal19}. However, the recent literature is mainly concerned with primary production
only and does not investigate the reasons why fault reactivation can occur during UGS
phases. Moreover, a very simplified geological structure is assumed in such analyses, with
a single fault in a two-dimensional (2D) vertical plane, and most likely this can only partially capture the
complex response expected from many intersecting faults in a fully three-dimensional (3D) environment
\cite{Zbi_etal17}. Only a few relatively old publications addressed the topic in UGS
reservoirs \cite{NagRoe97, Orl_etal13}. \citet{NagRoe97} developed a 2D geomechanical
model using the FLAC simulator for a typical faulted vertical section and concluded
that \textit{while the gas field is depleted, fault slip occurs due to compaction of the
reservoir and due to the upward movement of strata underlying the reservoir. Negligible
amounts of additional slip are induced when the reservoir is subjected to alternating
injection/extraction periods}. \citet{Orl_etal13} simulated the geomechanical behavior of
a specific UGS reservoir using the finite-element package DIANA. Their results highlighted
that \textit{the critically stressed section of the central fault affected by the fault
slipped~... during gas production. Additional fault slip could be expected during the
subsequent phase of cushion gas injection~... During annual cycles of gas injection and
production, the central fault is not critically stressed anymore}.

The aim of this work is multifold: i) to develop a robust computational framework 
allowing for the simulation of the inception of fault activation in
3D real-world geological settings; ii) to improve the understanding of the physical
mechanisms underlying induced seismicity during UGS activities, with specific reference to
the typical configurations of Dutch UGS reservoirs; iii) to investigate the factors that
can increase the chance of fault reactivation during UGS activities, identifying the
settings, conditions and material properties that could most likely cause ``unexpected''
fault reactivation in the reservoirs located in the Rotliegend formations;
iv) to define a set of practical guidelines allowing for a safe operational bandwidth in
such UGS fields, in consideration also of the different potential storage activities
(CH$_4$, H$_2$, CO$_2$). Some preliminary results have already been reported in the
scientific literature \cite{Tea_etal19, Tea_etal20}.

In order to accomplish such a complex multidisciplinary task, the overall work is
subdivided into two parts. The present paper (Part I) is mainly concerned with objectives
i) and ii), focusing on the mathematical and computational aspects of the modeling
approach, its application in a representative 3D test case of the problem of interest,
and the mechanisms that can cause ``unexpected'' fault activation during UGS
activities. 
\pt{The novelty of this contribution primarily lies in the identification and selection of a sophisticated state-of-the-art tool,
whose theoretical formulation is built on top of the works
\cite{franceschini2016novel, franceschini2019block, Iso_etal19}, capable of being effectively applied
to a realistic 3D reservoir case in the Rotliegend formation, compartmentalized by a complex fault system.
The application of this modeling framework provides us with the opportunity to explain, for the first time, the main mechanisms responsible for fault reactivation during CGI and UGS}.
A detailed sensitivity analysis for the different
storage activities and the definition of preliminary guidelines (aforementioned objectives iii) and iv)), is the target of Part II \cite{Tea_etal23}.

\oneway{In both Part I and Part II of this work, we adopt a one-way coupling strategy between fluid flow and geomechanics, which is acceptable for the conditions considered, given the weak interaction at the relevant scales. While this assumption simplifies the modeling and aligns with our focus on fault contact mechanics, we recognize its limitations and aim to extend the approach toward iterative or fully coupled schemes in future developments.}

The paper is organized as follows. The mathematical model of frictional contact mechanics
and flow in a 3D visco-elasto-plastic porous medium
is briefly introduced along with its numerical discretization and solution algorithms.
Faults are explicitly simulated within the porous rock as inner contact boundaries,
whose activation is macroscopically governed by Coulomb's criterion. Pressure change
within the faults, variation of Coulomb's parameters due to slip-weakening, and the
rheology of the caprock are properly accounted for.
The model is applied to a synthetic reservoir and fault system that realistically
represents the main geological features of the Rotliegend reservoirs.
%
Two scenarios are simulated to deepen the understanding of the geomechanical behavior of a
faulted UGS system. Computational results are presented and the mechanisms responsible for
fault reactivation during the UGS phases are outlined. A few conclusive remarks close the
publication.

\section{Methods and materials}
In this section, we discuss the development 
of the mathematical and numerical model used to
investigate the fault activation in the context of UGS
reservoirs. 
The aim is at solving the frictional contact mechanics problem for a faulted porous
medium, where the constraints are imposed in an exact way by Lagrange multipliers. The
friction behavior of the fracture is governed by Coulomb's criterion, with a
slip-weakening constitutive law. The variational formulation, its numerical discretization
and the possible related instability phenomena are discussed. The pore pressure, both in
the continuous matrix and inside the fracture network, is computed by a flow simulator
with a one-way coupled approach \cite{gambolati2000importance, Bui_etal19}, which turns
out to be fully warranted at the space and time scale of interest.
We use the quasi-static assumption, i.e., no acceleration contribution is accounted for,
under the hypothesis of likely negligible inertia of the system when small (e.g.,
centimetric) slip and small areal extent characterize the fault reactivation \cite{Biz10}.

\subsection{Strong formulation for the contact problem}
\begin{figure}
  \centering
  \includegraphics[width=0.25\textwidth]{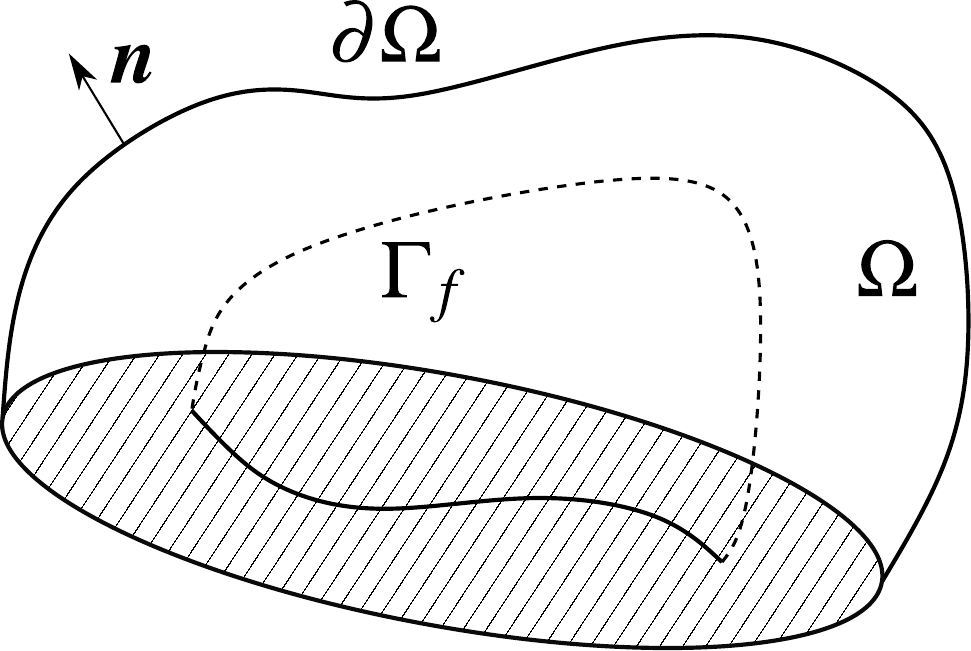}\qquad
  \includegraphics[width=0.25\textwidth]{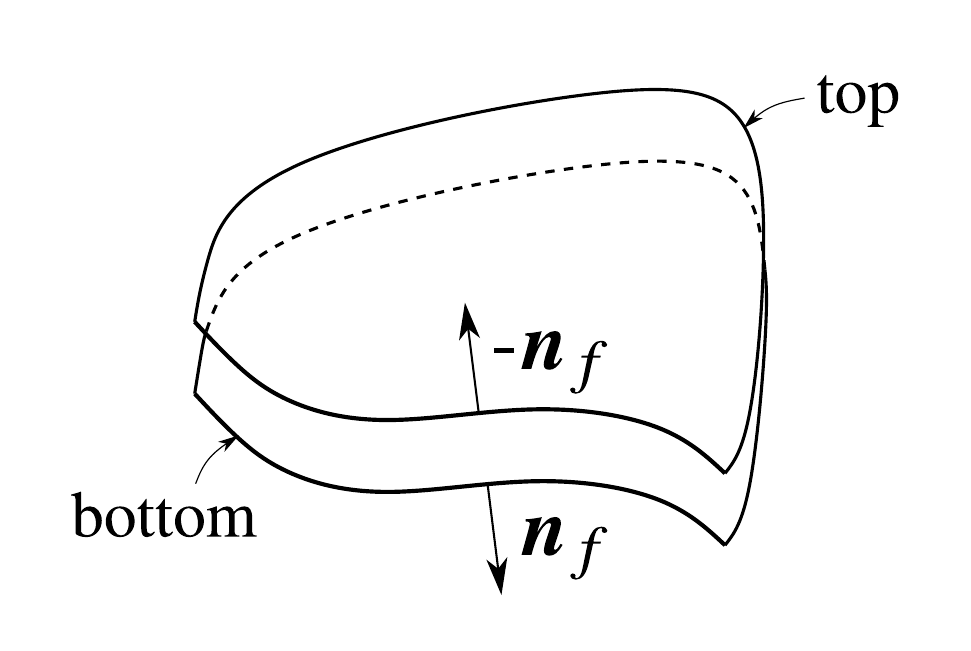}
  \caption{Sketch of the 3D domain $\Omega$ with its boundary, outer normal and inner fracture $\Gamma_f$ (left), made of the top and bottom contact surfaces and the normal direction $\vec{n}_f$ (right).}
  \label{fig:domain}
\end{figure}

A fault can be modeled at the macroscale as a lower dimensional internal boundary $\Gamma_f$ embedded in a 3D domain $\Omega \subset \mathbb{R}^3$. The fracture is represented as a
pair of surfaces in contact, conventionally denoted as \textit{top} and
\textit{bottom} and represented by $\Gamma_f^+$ and $\Gamma_f^-$, respectively. On such surfaces, normal and
frictional contact conditions have to be enforced, like the impenetrability of solid bodies and the
fulfillment of a friction criterion. In this work, we use Coulomb's frictional criterion
to provide the limiting modulus for the shear component of traction on the fault surface.
To complete the problem setting definition, we introduce the
external domain boundary $\Gamma \equiv \partial \Omega$, with its outer unit normal vector $\vec{n}$,
while $\vec{n}_f = \vec{n}_f^- = -\vec{n}_f^+$ denotes the normal direction to the fracture surface
$\Gamma_f$. Fig.~\ref{fig:domain} shows a sketch of the domain $\Omega$, the fault $\Gamma_f$ and the related quantities. Any
vector field can be decomposed along the normal and tangential direction to the fracture, i.e.,
$\vec{v} = v_N \vec{n}_f + \vec{v}_T$, with $v_N = \vec{n}_f^T \vec{v}$ and $\vec{v}_T =
\vec{v} - v_N \vec{n}_f = (\tensorTwo{1} - \vec{n}_f \otimes \vec{n}_f)\,\vec{v}$, where the
subscripts $N$ and $T$ are used to denote the normal and tangential components,
respectively, and $\tensorTwo{1}$ is the identity tensor of order 2.

Assuming quasi-static conditions and
infinitesimal strains, the strong form of the linear momentum balance at every instant $t$ in the time interval $[0,t_{\max}]$ can be stated as follows
\cite{kikuchi1988contact, laursen2003computational, wriggers2006computational}: find the
displacement vector ${\vec{u}: \overline{\Omega}\times[0,t_{\max}] \rightarrow \mathbb{R}^3}$ such that:
\begin{subequations}
\begin{align}
  \nabla \cdot \hat{\tensorTwo{\sigma}} (\vec{u}) + \vec{b} &= 0 & &\text{in } \Omega \times [0,t_{\max}],
    \label{eq:mom_bal}\\
  \vec{u} &= \bar{\vec{u}} & &\text{on } \Gamma_{u} \times [0,t_{\max}], \\
  \hat{\tensorTwo{\sigma}} (\vec{u}) \cdot \vec{n} &= \bar{\vec{t}} & &\text{on }
    \Gamma_{\sigma} \times [0,t_{\max}],
\end{align}
\label{eq:governing_eqs}\null
\end{subequations}
where $\hat{\tensorTwo{\sigma}}$ is the total stress tensor, $\vec{b}$ collects the
external body loads and $\Gamma_{u} \cup \Gamma_{\sigma} = \Gamma$, $\Gamma_{u} \cap
\Gamma_{\sigma} = \varnothing$, are the portion of the boundary where Dirichlet and
Neumann conditions are imposed, respectively. On the fracture $\Gamma_f$, normal and
friction compatibility conditions need to be enforced \cite{laursen2003computational,
wriggers2006computational}. The normal contact conditions on the fracture read:
\begin{subequations}
\begin{align}
  t_N = \vec{t} \cdot \vec{n}_f &\le 0 & &\text{only compressive traction is allowed}, \\
  g_N = \llbracket \vec{u} \rrbracket \cdot \vec{n}_f &\ge 0 & &\text{impenetrability
    condition}, \\
  t_N g_N &= 0 & &\text{either the fracture is compressed or it is open}.
\end{align}
\label{eq:contact_cond_N}\null
\end{subequations}
The conditions for the frictional component are:
\begin{subequations}
\begin{align}
  \left\| \vec{t}_T \right\|_2 \le \tau_{\max}\left(t_N,\|\vec{g}_T\|_2\right)& &
    &\text{Coulomb's criterion}, \\
  \dot{\vec{g}}_T \cdot \vec{t}_T - \tau_{\max}\left(t_N,\|\vec{g}_T\|_2\right)
    \|\dot{\vec{g}}_T\|_2 = 0& & &\text{frictional traction is aligned with sliding rate}
    \label{eq:contact_cond_T2}.
\end{align}
\label{eq:contact_cond_T}\null
\end{subequations}
In Eqs.~\eqref{eq:contact_cond_N}-\eqref{eq:contact_cond_T}, we split the traction $\vec{t}$ on the
fracture and the displacement jump across it into normal and tangential components, i.e.,
$\vec{t} = t_N \vec{n}_f + \vec{t}_T$ and $\llbracket \vec{u} \rrbracket = g_N \vec{n}_f +
\vec{g}_T$, respectively. The jump is defined as $\llbracket \vec{u} \rrbracket =
\vec{u}|_{\text{top}} - \vec{u}|_{\text{bottom}}$. To characterize the standard Coulomb frictional
criterion, $c$ and $\varphi$ are introduced, i.e., the cohesion and the friction angle,
respectively, obtaining:
\begin{equation}
  \tau_{\max}(t_N,\|\vec{g}_T\|_2) = c - t_N\tan\left( \varphi\left(
    \|\vec{g}_T\|_2\right)\right).
  \label{eq:Coulomb}
\end{equation}
In Eq.~\eqref{eq:Coulomb}, the friction angle generally depends on the modulus of the tangential
component of the displacement jump, i.e., the slippage, so as to simulate a slip-weakening frictional behavior.
Since a quasi-static approach is used, we can 
replace the tangential
displacement rate $\dot{\vec{g}}_T$ in Eq.~\eqref{eq:contact_cond_T2} with the incremental
tangential displacement $\Delta\vec{g}_T$ with respect to the previous time-step value \cite{wohlmuth2011variationally}.

The fault surface $\Gamma_f$ can be split into three non-intersecting portions. Each portion is characterized by a different operating mode allowed for by the possible
combinations of the previous conditions~\eqref{eq:contact_cond_N}-\eqref{eq:contact_cond_T}:
\begin{itemize}
  \itemsep 0em
  \item \textit{Stick}: the surface is compressed ($t_N < 0$) and the shear
traction modulus does not exceed the limiting value provided by the criterion in Eq.~\eqref{eq:Coulomb}. When there has been no prior sliding, the displacement field remains continuous across $\Gamma_f$. However, if sliding has occurred previously, there will be no additional tangential slip increment;
  \item \textit{Slip}: the normal traction is still negative, but the surface is free to
slip. In this case, Coulomb's equality holds \cite{simo1998computational}:
\begin{equation}
  \vec{t}_T^* = \tau_{\max}(t_N,\|\vec{g}_T\|_2)
\frac{\Delta\vec{g}_T}{\|\Delta\vec{g}_T\|_2}.
\end{equation}
  Only the normal component of the displacement field is continuous across $\Gamma_f$;
  \item \textit{Open}: the normal traction is non-negative and the two contact surfaces $\Gamma_f^+$ and $\Gamma_f^-$ are free to
move on condition to avoid compenetration. Hence, the displacement field across $\Gamma_f$ is discontinuous and the traction on the fault
vanishes, i.e., $\vec{t} = 0$.
\end{itemize}
For additional details on the mathematical formulation, see \cite{kikuchi1988contact,
laursen2003computational, wriggers2006computational} and more recently
\cite{franceschini2016novel, franceschini2020algebraically, franceschini2022scalable}.

According to Terzaghi's principle, the total stress tensor in a saturated porous medium
can be decomposed as the sum of two contributions: the effective stress tensor acting on the solid
skeleton and a volumetric term depending on the averaged fluid pressure $p$:
\begin{equation}
  \hat{\tensorTwo{\sigma}} = \begin{cases}
    \tensorTwo{\sigma} - \tensorTwo{1} p & \text{on } \Gamma_{\sigma} \cup \Gamma_f,\\
    \tensorTwo{\sigma} - \alpha \tensorTwo{1} p = \tensorFour{C}  :
    \tensorTwo{\varepsilon} - \alpha \tensorTwo{1} p & \text{in } \Omega,
  \end{cases}
  \label{eq:terz_princ}
\end{equation}
where the fluid pressure is averaged by the saturation indices of the different phases, $\alpha$ is the Biot coefficient, taking care of the ratio
between the grain and the porous matrix compressibility, and $\tensorTwo{1}$ the identity tensor of order $2$
\cite{coussy2004poromechanics}. The effective stress tensor $\tensorTwo{\sigma}$ is called effective Terzaghi stress tensor on the domain boundaries, while it is the effective Biot stress tensor inside the domain itself. In Eq.~\eqref{eq:terz_princ}, the constitutive relationship
defining the effective Biot stress tensor is introduced, where $\tensorFour{C}$ is a
fourth order elasticity tensor, generally non-linear, and $\tensorTwo{\varepsilon} = \nabla^s \vec{u}$ is the
strain tensor, with $\nabla^s = (\nabla + \nabla^T)/2$ the symmetric gradient operator. The
mechanical constitutive law relates a strain variation in the porous medium to an effective stress variation. Such a law can
be described by a simple linear elastic model (Hooke's law), with constant or variable parameters, but
also more complex elasto-plastic rules with time-dependent contributions can be introduced, e.g., a visco-elasto-plastic law. For
more details on the appropriate constitutive laws and their implementations, see
\cite{simo1998computational, de2011computational}.

\subsection{Mass balance equation}
The mass conservation of the fluid species $\kappa$ reads \cite{xikui1992multiphase,
rutqvist2002modeling, yadigaroglu2018introduction, wilkins2021open}:
\begin{equation}
  \frac{\partial}{\partial t}(\rho^\kappa) + \nabla \cdot \vec{F}^\kappa = q_s^\kappa,
  \label{eq:mass_bal}
\end{equation}
where $\rho^\kappa$ and $\vec{F}^\kappa$ are the density and the flux, respectively, of the fluid species
$\kappa$. The density $\rho^\kappa$ represents the mass of $\kappa$ per unit of rock volume and can be written as:
\begin{equation}
  \rho^\kappa = \phi \sum_{\beta} S_\beta \rho_\beta \chi_\beta^\kappa,
  \label{eq:dens_k}
\end{equation}
with $\phi$ the porosity, $S_\beta$ the saturation of phase $\beta$,
that can be either liquid or gas, $\rho_\beta$ the density of phase $\beta$, 
and $\chi_\beta^\kappa$ the mass
fraction of component $\kappa$ in phase $\beta$. In this application, we assume isothermal conditions, meaning that the fluid density is primarily influenced by pressure. However, it can also be affected by other factors, such as mass fraction, according to a specific equation of state. Saturations and mass fractions are constrained by the well-known conditions:
\begin{equation}
  \sum_{\beta} S_\beta = 1 \quad \text{and} \quad \sum_\kappa \chi_\beta^\kappa = 1 \quad
  \forall \beta.
  \label{eq:constrSx}
\end{equation}
The fluid flux of component $\kappa$ is the sum of the fluxes for each phase:
\begin{equation}
  \vec{F}^\kappa = \sum_\beta \chi_\beta^\kappa \vec{F}_\beta,
  \label{eq:Fkappa}
\end{equation}
and each phase flux is described by Darcy's law as:
\begin{equation}
  \vec{F}_\beta= \rho_\beta \vec{v}_\beta = -\rho_\beta \frac{\tensorTwo{k}\
k_{r,\beta}}{\mu_\beta}\left(\nabla p_\beta - \rho_\beta \vec{g}\right),
  \label{eq:Fbeta}
\end{equation}
with $\mu_f$ and $k_{r,\beta}$ the viscosity and the relative permeability of the phase,
$\tensorTwo{k}$ the permeability tensor, and $p_\beta$ the pressure in phase $\beta$.

Among the various models available that consider poro-elastic effects for updating
porosity, e.g., see \cite{kim2012formulation, martinez2013coupled}, in this study we
follow the approach presented in these works \cite{coussy2004poromechanics,
garipov2016discrete}, which expresses porosity as:
\begin{equation}
  \phi = \phi_0 + \alpha \varepsilon_v +
\frac{\left(\alpha-\phi_0\right)\left(1-\alpha\right)}{K_d}\left(p-p_0\right),
  \label{eq:phidef}
\end{equation}
where $K_d$ is the drained bulk modulus, 
$\varepsilon_v = \text{trace}\left(\tensorTwo{\varepsilon}\right)$ is the
volumetric strain, and $\phi_0$ and $p_0$ are the reference porosity and fluid pressure
respectively. In Eq.~\eqref{eq:phidef}, $p$ is the averaged fluid pressure, computed as:
\begin{equation}
  p = \sum_\beta S_\beta p_\beta.
  \label{eq:pfluid}
\end{equation}

In oedometric conditions and a constant total stress state, from Terzaghi's
principle (Eq.~\eqref{eq:terz_princ}) we have $d\sigma_z = \alpha\, dp$ and
\begin{equation}
  \frac{\partial}{\partial t}\varepsilon_v = 
    C_m\alpha\frac{\partial}{\partial t}p,
  \label{eq:depsv}
\end{equation}
with $C_m$ is the vertical uniaxial
compressibility. 
In this case,
the mass balance
is decoupled from the linear momentum balance and can be solved in advance, providing
a pressure field acting as an external body load for the structural problem.

%
Even though this assumption is not strictly guaranteed for the reservoir application considered in the present work, at the (large) spatial and temporal scales of interest, the coupling between flow and mechanics is weak. Thus, a one-way coupled approach, where Eqs.\eqref{eq:mass_bal} are first solved for all fluid species and the averaged pressure from Eq.\eqref{eq:pfluid} is subsequently introduced into Eq.~\eqref{eq:governing_eqs}, is fully justified; see, for instance, \cite{gambolati2000importance, zoccarato2016data, teatini20143d, castelletto2013geological, castelletto2013compartmentalization, castelletto2013multiphysics, rahman2022effect}.
\oneway{To the best of the authors' knowledge, few studies adopt a fully coupled simulation framework, and those that do are typically limited to small to medium-scale cases \cite{mathur2024thermo} or rely on iterative approaches \cite{nasrollahzadeh2021field}.
Recent studies have also compared one-way and two-way coupling schemes under conditions similar to those of this work, further supporting our approach \cite{an2025study, ferrari2025effects}. Therefore, a one-way coupling is adopted here. This choice is additionally motivated by the fact that, at reservoir depth, the fluids involved are generally more compressible than the host rock, and our primary research objective is to investigate fault contact mechanics rather than the full complexity of flow-deformation interactions. Nonetheless, we acknowledge the limitations of the one-way approach and aim to extend this work toward an iterative coupling scheme in future developments.}

The set of Eqs.~\eqref{eq:mass_bal} for obtaining the pressure field in the porous medium is usually solved applying a finite volume
method because it preserves the mass conservation at the elemental level~\cite{eymard2000finite, manzini2008finite, droniou2014finite}. Nevertheless, also a finite
element or mixed finite element approach can be successfully used, e.g.,~\cite{wan2003stabilized, jha2007locally, prevost2016faults, Abushaikha2017, NarFerAbu23}.
In our analysis, the numerical simulation of the
multiphase flow in the porous matrix has been carried out by using \textit{OPM Flow}, an open-source
reservoir simulator based on a classical finite volume discretization~\cite{OPM, rasmussen2021open}.
As to the computation of the pressure field within the network of faults, two strategies can be employed: either the
domain explicitly contains the faults as \textit{thin} 3D cells, or the pressure is extended to the faults from the surrounding 3D
cells, according to some physical treatment of the contact surfaces as inner boundaries. We decided to use the latter approach, thus allowing to represent the fault at the macroscale as zero-thickness lower dimensional elements. Generally speaking, the two limiting cases that can be met in reality are
{\em sealing} and {\em non-sealing} faults. In the former case, the fault acts as an impermeable barrier and the pressure change does not propagate
from one side to the other of the contact surfaces. In this situation, we can assume that the pressure variation in the fault is null. 
In contrast, in the latter case, the fault is fully permeable and does not exhibit any resistance to fluid flow. In this situation, we assume the
pressure variation in the fault to be equal to the arithmetic average of the pressure computed on the two side cells. 

\subsection{Variational formulation and discretization}
In this section, the variational formulation for the strong form of the linear momentum
balance in Eq.~\eqref{eq:governing_eqs}, equipped with the constraints of
Eqs.~\eqref{eq:contact_cond_N}-\eqref{eq:contact_cond_T}, is presented.
The weak form of the governing equations naturally produces a variational inequality because of the frictional contact constraints~\cite{kikuchi1988contact}. In order to avoid this difficulty, it is possible to reduce
the original inequality to a standard variational formulation by an active-set strategy and
either a penalty regularization or the introduction of Lagrange multipliers. We
decided to use the Lagrange multiplier technique, which can be computationally more expensive, because new primary
unknowns are introduced and the resulting algebraic problem gains a saddle-point nature, but generally much more accurate, stable and robust, not requiring additional parameters. Moreover, though generating saddle-point systems, this formulation allows to produce a sequence of linear problems less sensible to ill-conditioning
issues \cite{ferronato2012parallel}. From a physical viewpoint, Lagrange
multipliers represent the traction field on the fault surfaces, thus the stress evaluation
and the related reactivation hazard becomes straightforward.

\af{We emphasize that to retrieve the complete variational formulation for the problem at
hand, discussing the well-posedness of the problem and defining the right functional spaces and conditions, is beyond the purposes of
this work. Here, we use the main theoretical outcomes obtained in
\citet{kikuchi1988contact} and \citet{wohlmuth2011variationally}, further discussed in
other works \cite{wriggers2006computational, franceschini2020algebraically}. 
The selected function spaces and the residual equations that guarantee the well-posedness of the problem are reported below.} The notations $(\cdot,\cdot)_\Omega$ and
$\langle\cdot,\cdot\rangle_\Gamma$ denote the $L^2$-inner product of functions in $\Omega$
(3D domain) and $\Gamma$ (lower dimensional domain), respectively. Let $\mathcal{V} =
[H^1(\Omega)]^3$ be the Sobolev space of vector functions whose first derivatives belong
to $L^2(\Omega)$; let $\mathcal{M}$ be the dual space of the trace space $\mathcal{W} =
[H^{1/2}(\Omega)]^3$; and let $\mathcal{M}(t_N,\|\vec{g}_T\|)$ be its subspace such that
\begin{equation}
  \mathcal{M}(t_N,\|\vec{g}_T\|) = \left\{\vec{\mu} \in \mathcal{M}: \langle
\vec{\mu},\vec{v} \rangle_{\Gamma_f} \le \langle \tau_{\max}(t_N,\|\vec{g}_T\|_2),
\|\vec{v}_T\|\rangle_{\Gamma_f}, \vec{v} \in \mathcal{W} \text{ with } v_N \le 0\right\}.
  \label{eq:Mdef}
\end{equation}

Given the finite-dimensional subspaces $\mathcal{V}^h \subset \mathcal{V}$ and
$\mathcal{M}^h(t^h_N,\|\vec{g}^h_T\|) \subset \mathcal{M}(t_N,\|\vec{g}_T\|)$, the finite dimensional weak form
of the problem in Eq.~\eqref{eq:governing_eqs} with Terzaghi relation of Eq.~\eqref{eq:terz_princ} and conditions
Eqs.~\eqref{eq:contact_cond_N}-\eqref{eq:contact_cond_T} can be stated as follows: at every instant $t\in[0,t_{\max}]$, find
$\{\vec{u}^h,\vec{t}^h\} \in \mathcal{V}^h \times
\mathcal{M}^h(t^h_N,\|\vec{g}^h_T\|)$ such that:
\begin{subequations}
\begin{align}
  \mathcal{R}_u &= \left(\nabla^s \vec{\eta},\hat{\tensorTwo{\sigma}}\right)_{\Omega}
- \langle \vec{\eta},\hat{\tensorTwo{\sigma}} \cdot \vec{n}\rangle_{\Gamma}
- \left(\vec{\eta},\vec{b}\right)_{\Omega} \nonumber \\
  &= \left(\nabla^s \vec{\eta},\hat{\tensorTwo{\sigma}}\right)_{\Omega}
- \langle  \vec{\eta},\hat{\tensorTwo{\sigma}} \cdot \vec{n}_f^+\rangle_{\Gamma_f^+}
- \langle  \vec{\eta},\hat{\tensorTwo{\sigma}} \cdot \vec{n}_f^-\rangle_{\Gamma_f^-}
- \langle \vec{\eta}, \hat{\tensorTwo{\sigma}} \cdot \vec{n}\rangle_{\Gamma_{\sigma}}
- \left(\vec{\eta},\vec{b}\right)_{\Omega} \nonumber \\
  &= \left(\nabla^s \vec{\eta},\tensorTwo{\sigma}(\vec{u}^h) - \alpha \tensorTwo{1} p\right)_{\Omega}
- \langle \llbracket \vec{\eta} \rrbracket,\vec{t}^h - p \vec{n}_f\rangle_{\Gamma_f}
- \langle \vec{\eta}, \bar{\vec{t}}\rangle_{\Gamma_{\sigma}}
- \left(\vec{\eta},\vec{b}\right)_{\Omega} = 0 & &\forall \vec{\eta} \in \mathcal{V}^h,
\label{eq:ResU}\\
  \mathcal{R}_t &= \langle t_N^h - \mu, g_N \rangle_{\Gamma_f} + \langle
\vec{t}_T^h-\vec{\mu}_T, \Delta \vec{g}_T\rangle_{\Gamma_f} \ge 0 & &\forall \vec{\mu} \in
\mathcal{M}^h(t^h_N,\|\vec{g}^h_T\|). \label{eq:ResT}
\end{align}
\label{eq:ResFE}\null
\end{subequations}
We use a Galerkin approach, hence the test functions $\vec{\eta}$ and
$\vec{\mu}$ belong to the same
function spaces used to define the trial functions for the displacement and traction field, respectively. To transform the variational
inequality of Eq.~\eqref{eq:ResT} into a variational equality, an iterative
active-set algorithm \cite{nocedal2006numerical, antil2018frontiers} is applied.
According to this approach, the fault surfaces $\Gamma_f$ is subdivided into active and inactive
regions for both components of the traction, i.e., the {\em Stick} ($\Gamma_f^{\text{stick}}$, active for normal and tangential components), {\em Slip} ($\Gamma_f^{\text{slip}}$, active for normal component), and {\em Open} ($\Gamma_f^{\text{open}}$, inactive) portions of $\Gamma_f$. With this subdivision, the
variational inequality of Eq.~\eqref{eq:ResT} becomes:
\begin{align}
  \mathcal{R}_t = \langle \vec{\mu}, \llbracket\vec{u}^h\rrbracket \rangle_{\Gamma^{\text{stick}}_f} + \langle
\mu_N, g_N \rangle_{\Gamma^{\text{slip}}_f} + \frac{1}{k}\langle \vec{\mu}_T, \vec{t}_T^h
- \vec{t}_T^* \rangle_{\Gamma^{\text{slip}}_f} + \frac{1}{k} \langle \vec{\mu}, \vec{t}^h
\rangle_{\Gamma^{\text{open}}_f} &= 0 & &\forall \vec{\mu} \in
\mathcal{M}^h(t^h_N,\|\vec{g}^h_T\|),
  \label{eq:ResTeq}
\end{align}
where $k$ is a unitary coefficient introduced to ensure dimensional consistency.
The non-linear system of Eqs.~\eqref{eq:ResU}-\eqref{eq:ResTeq} is solved by a Newton linearization and at convergence we check the consistency of the traction state on the faults with the initial subdivision of $\Gamma_f$ into $\Gamma_f^{\text{stick}}$, $\Gamma_f^{\text{slip}}$, and $\Gamma_f^{\text{open}}$. If the consistency check is satisfied, the active-set algorithm is stopped and we can move to the following time instant, otherwise a new subdivision of $\Gamma_f$ is defined, Eq.~\eqref{eq:ResTeq} is re-computed and the resulting non-linear system solved again.

Introducing $\vec{u}^h = \sum_i u_i \vec{\eta}_i$ and $\vec{t}^h = \sum_j t_j
\vec{\mu}_j$, i.e., the discrete representation of the displacement and traction fields,
where $\{\vec{\eta}_i\}$ and $\{\vec{\mu}_j\}$ are bases for $\mathcal{V}^h$ and
$\mathcal{M}^h(t^h_N,\|\vec{g}^h_T\|)$, respectively, the set of variational equalities of Eqs.~\eqref{eq:ResU}-\eqref{eq:ResTeq}
becomes an algebraic nonlinear system. 
The bases for the finite-dimensional spaces $\mathcal{V}^h$ and $\mathcal{M}^h$ are selected with the aid of the finite element method.
Given the regularity requirements defined above, we use low-order discretization spaces for both displacement and traction.
The computational domain is subdivided into non-overlapping hexahedral
elements, $\Omega = \bigcup_{i=1}^{n_e}\Omega_i$ and $\Omega_i\cap\Omega_j=\varnothing$ for any $i\neq j$. This choice is done for the sake of the consistency with the domain discretization used for the multiphase flow model, which is based on a standard finite volume approach.
We use a conformal representation of
the faults, i.e., $\Gamma_f = \bigcup_j\partial \Omega_j$, where $\Omega_j$ are elements
sharing a face with $\Gamma_f$. In such a way, the fault contact surfaces are composed by pairs of quadrilateral
elements. Each pair of quadrilateral elements is also denoted as a zero-thickness {\em interface finite element}~\cite{goodman1968model}. 
According to the value of the traction, every interface element can change its status, i.e., it can belong to either the
stick, slip or open portion of $\Gamma_f$.

The mixed finite element discretization adopted in this work is described in these works
\cite{franceschini2020algebraically, franceschini2022scalable}. It consists of a
$\mathbb{Q}_1$ first-order interpolation for the nodal-based displacement field and a
$\mathbb{P}_0$ piecewise constant interpolation for the element-based traction field. 
The collection of coefficients $u_i$ and $t_j$ are
the components of the unknown algebraic arrays, named $\Vec{u}$ and $\Vec{t}$ of sizes $3 n_n$ and $3 n_f$, with $n_n$
the number of nodes in the hexahedral mesh and $n_f$ the number of interface elements.
This
approach has the main advantage of being naturally coupled with a finite volume pressure
solution computed on the same domain discretization at no additional cost, since both the traction and the
pressure are represented using the same space on the same grid. On the other side, in order to
ensure the LBB-stability of the proposed mixed finite element spaces, a tailored jump
stabilization has been proposed in Franceschini et al. \cite{franceschini2022scalable}. For a complete
analysis of the LBB-stability of a general pair of mixed finite element spaces, see
\citet{elman2014finite}. An implementation of the presented algorithm can be found in this
website \cite{ATLAS}.

%
We emphasize that even if a linear elastic constitutive relation is used for the porous
medium, the set of equations reported in Eqs.~\eqref{eq:ResU}-\eqref{eq:ResTeq} represents a non linear problem,
because a consistent partitioning of the fracture surface is unknown and has to be
computed depending on the solution vectors. To be more specific, constraints in
Eqs.~\eqref{eq:contact_cond_N}-\eqref{eq:contact_cond_T} are the Karush-Kuhn-Tucker (KKT)
conditions and we are dealing with a non-linear optimization problem
\cite{nocedal2006numerical}. As already mentioned, the solution strategy is based on an
active-set strategy, with each non-linear problem addressed by a classical exact Newton algorithm.


\subsection{Linearization and linear system solution}

The use of a mixed finite element approximation produces a linearized
step with a generalized saddle-point Jacobian matrix~\cite{benzi2005numerical}.
The Jacobian is generally non-symmetric because of the contribution related to the friction component of the
traction when the fracture slides.
In particular, at each Newton iteration, the linear system that has to be solved is:
\begin{equation}
  J \delta \Vec{x} = -\Vec{r},
\end{equation}
with the $2\times 2$ block matrix $J$, the residual vector $\Vec{r}$, and the solution
vector $\delta\Vec{x}$ given by:
\begin{equation}
  J = \begin{bmatrix}
    \frac{\partial \mathcal{R}_u}{\partial \vec{u}} & \frac{\partial
\mathcal{R}_u}{\partial \vec{t}} \\
    \frac{\partial \mathcal{R}_t}{\partial \vec{u}} & \frac{\partial
\mathcal{R}_t}{\partial \vec{t}}
  \end{bmatrix}, \quad \Vec{r} = \begin{bmatrix}
    \mathcal{R}_u \\
    \mathcal{R}_t
  \end{bmatrix}, \quad \text{and} \quad \delta\Vec{x} = \begin{bmatrix}
    \delta\Vec{u} \\
    \delta\Vec{t}
  \end{bmatrix},
\end{equation}
where $\mathcal{R}_u$ and $\mathcal{R}_t$ are computed at the current counter level $l$.
As usual, the updated solution vector at level $l+1$ is:
\begin{equation}
  \begin{bmatrix}
    \Vec{u} \\
    \Vec{t} \\
  \end{bmatrix}^{l+1} = \begin{bmatrix}
    \Vec{u} \\
    \Vec{t} \\
  \end{bmatrix}^l + \begin{bmatrix}
    \delta\Vec{u} \\
    \delta\Vec{t} \\
  \end{bmatrix}.
\end{equation}
For the detailed expression of the Jacobian, see \cite{franceschini2020algebraically}. The
resulting linear system is characterized by a large and sparse matrix and it is necessary
to use a preconditioned iterative method for its efficient solution. To achieve satisfactory results, the use of a
suitable preconditioner is mandatory. Since the properties of the linear system may change
significantly as the simulation proceeds and the fault elements change status, also the
preconditioner must evolve. Among others, an idea is to exploit the scalability
intrinsically present in the multigrid approach and combine it with the known
physics-based partitioning of the blocks to be able to solve the saddle-point matrix. For
details on robust and efficient techniques used for the solution of this peculiar linear
system, the reader may refer to the specific literature \cite{franceschini2019block,
franceschini2022scalable, franceschini2022reverse}.

%
%
\subsection{Constitutive model for fracture: slip weakening}
\begin{figure}
  \centering
  \null\hfill
  \includegraphics[width=0.33\textwidth]{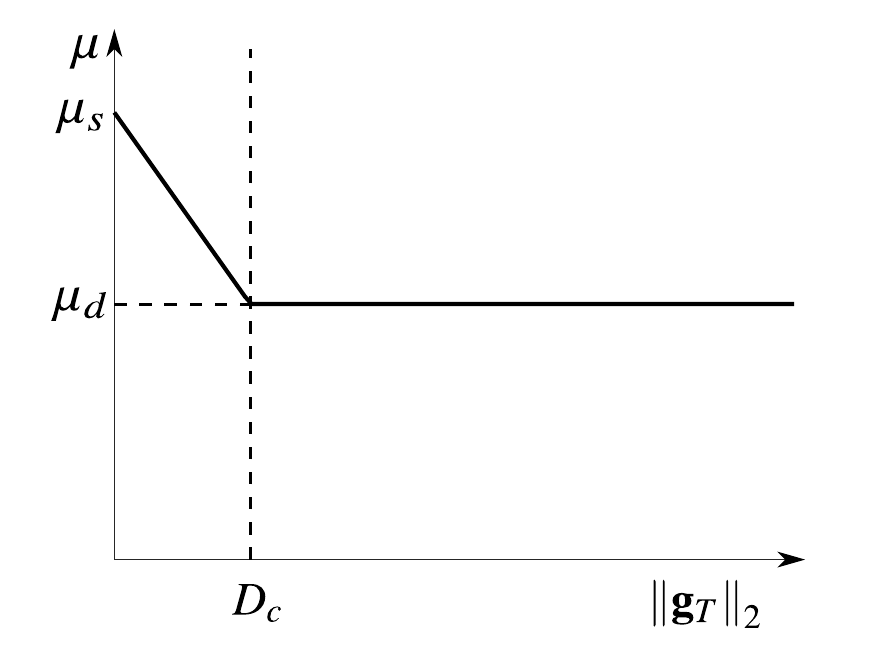}\hfill
  \includegraphics[width=0.33\textwidth]{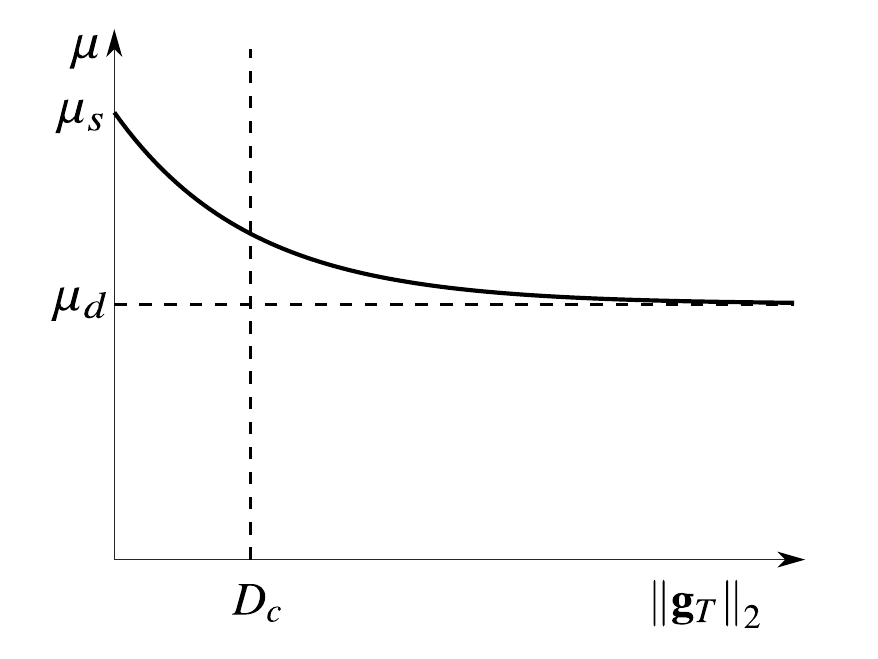}\hfill
  \includegraphics[width=0.33\textwidth]{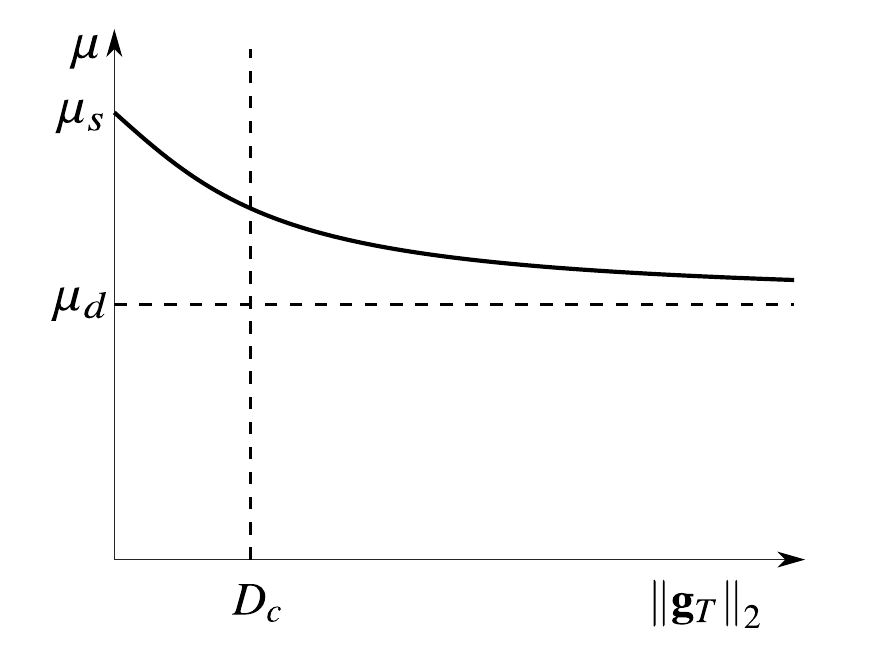}\hfill\null
  \caption{Slip weakening constitutive laws. From left to right: piecewise linear friction
law, exponential friction law and inverse trigonometric friction law.}
  \label{fig:weak_laws}
\end{figure}

In this work we use both the classical Coulomb criterion with a constant friction
coefficient, and a slip-weakening friction law with a variable
friction coefficient. Originally used by \citet{andrews1976rupture} to take into account the change from static to dynamic friction, slip-weakening friction laws \cite{ida1972cohesive,palmer1973growth}
 are based on the concept that the shear stiffness of
the fracture decreases as sliding occurs. From a mathematical viewpoint, the standard
Coulomb criterion reads:
\begin{equation}
  \left\|\vec{t}_T\right\|_2 \le c - t_N \mu, \quad \text{with } \mu = \tan\varphi,
\end{equation}
while a more general slip-weakening friction law reads:
\begin{equation}
  \left\|\vec{t}_T\right\|_2 \le c - t_N \mu(\|\vec{g}_T\|_2).
\end{equation}
A simple expression to account for the friction reduction with fault motion is provided by a piecewise linear function, as shown in
Fig.~\ref{fig:weak_laws}, where the friction coefficient linearly decreases from the static value $\mu_s$ down to the dynamic value $\mu_d$ at a
sliding value equal to $D_c$. For larger sliding values, the friction coefficient remains constantly equal to $\mu_d$.
Other analytical expression can be also used to simulate the friction
coefficient reduction with the sliding, such as an exponential law (see Fig.~\ref{fig:weak_laws}), which has the advantage to allow for a
\textit{smooth} variation that can be differentiated everywhere. 
It provides:
\begin{equation}
  \mu = \mu_d + \left(\mu_s - \mu_d\right)
    \exp\left(-\frac{\left\|\vec{g}_T\right\|_2}{D_c}\right).
\end{equation}
A similar smooth behavior can be formulated based on inverse trigonometric functions (see Fig.~\ref{fig:weak_laws}):
\begin{equation}
  \mu = \mu_d + \left(\mu_s - \mu_d\right)
    \left(1-\frac{2}{\pi}\arctan\frac{\left\|\vec{g}_T\right\|_2}{D_c}\right).
\end{equation}

\begin{figure}
  \centering
  \includegraphics[width=0.8\textwidth]{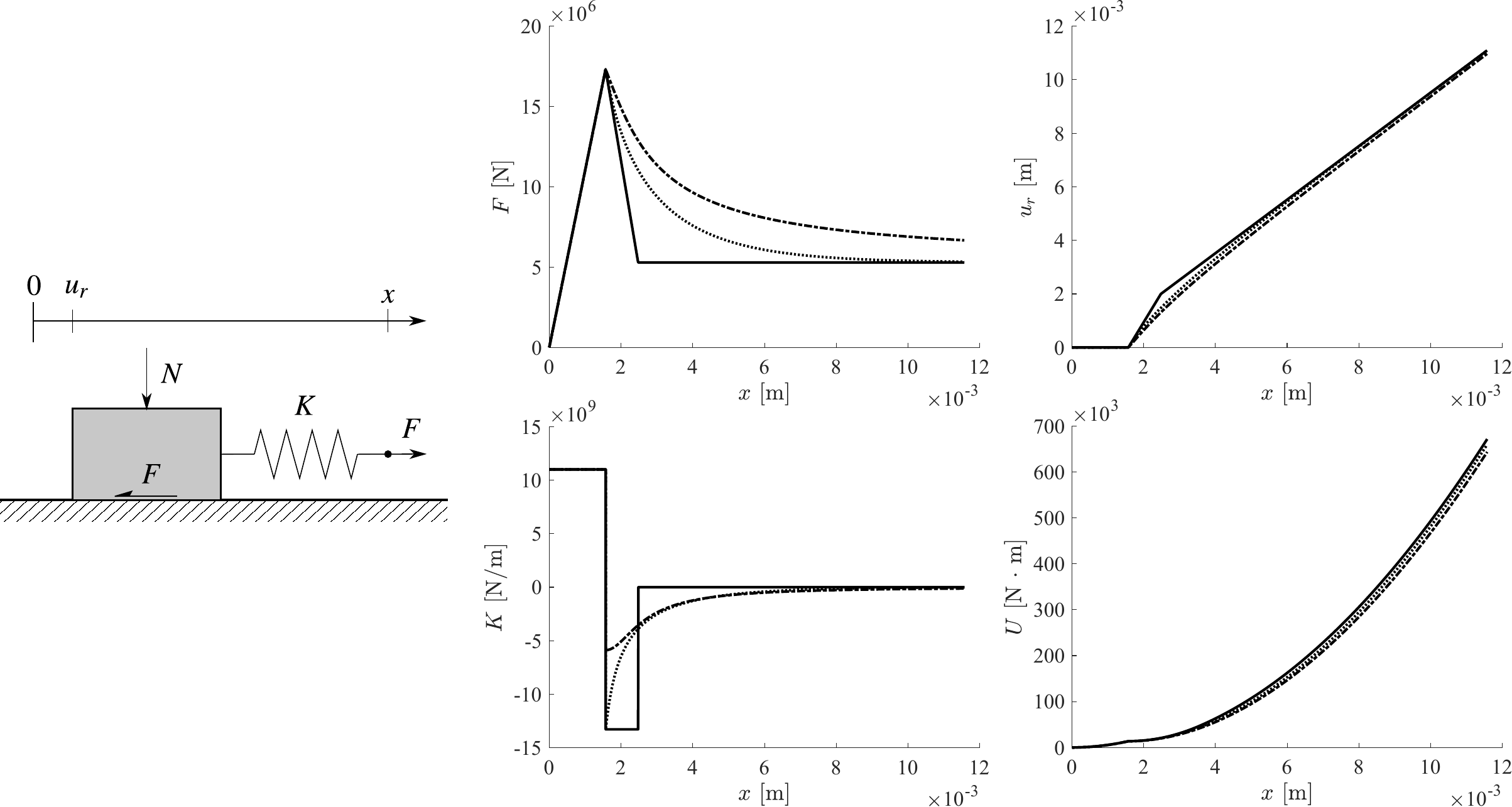}
  \caption{Mono-dimensional friction system with 1 degree of freedom. On the left: sketch
of the model used in this example. On the right, from the top left, (i) friction strength,
(ii) relative displacement, (iii) global system stiffness and (iv) internal energy.
Continuous, dotted and dashed line represent linear, exponential and inverse trigonometric
slip-weakening formulation, respectively.}
  \label{fig:fric_1D}
\end{figure}

In order to compare the three different expressions, we use the simple 1D problem sketched in
Fig.~\ref{fig:fric_1D}. The selected physical parameter set is: $\mu_s = \tan(30^\circ)$,
$\mu_d = \tan(10^\circ)$, $D_c = 2$~mm, with the spring stiffness $K = 11\times10^9$~N/m and a compression load $N = 3\times10^7$~N. The
first three values are representative of the conditions typically found in the seismogenic gas fields within the Rotliegend
stratigraphic units in the Netherlands~\cite{Bui_etal19,Hun_etal20}. The physical
quantities of interest are shown in Fig.~\ref{fig:fric_1D}. The primary variable is always
the displacement of the point where the external load $N$ is applied, while the outcomes are:
(i) the friction reaction $F$, (ii) the relative displacement $u_r$ between the body connected to the spring and the fixed basement, (iii) the global system
stiffness $\overline{K} = \partial F / \partial u_r$, and (iv) the internal energy $U = \int_0^{u_r} F\, dx$. Though the response in terms of friction
strength are different, both relative displacement and global energy are
comparable. By distinction, the global stiffness behaves differently and for two cases out of
three it reaches negative amounts that are greater in absolute value than the original spring
stiffness $K$. 
The finite element approach used in the present modeling analysis is based on the global
equilibrium of the system and not on a local (elemental) balance. This is the reason why a comparison of
the global energy associated with the different laws is meaningful. At the elemental level, it is desirable to avoid negative stiffness, 
which could potentially lead to
friction instabilities. Hence, we chose to work
with the law based on the inverse trigonometric function, i.e., the only providing a
minimum stiffness smaller in absolute value than the original one. In such a way,
we can ensure, at least for conditions similar to the ones used in this example, a
positive global stiffness. 

%

\section{Analyses and results}
\subsection{Conceptual model}
\begin{figure}
  \centering
  \null\hfill
  \includegraphics[height=0.3\textwidth]{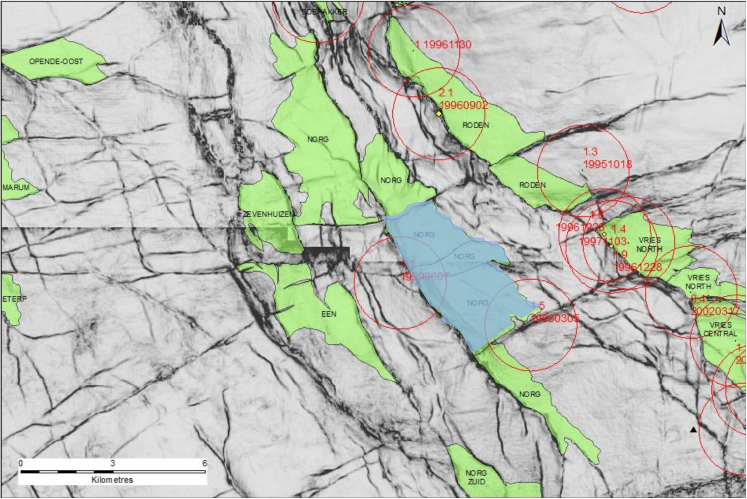}\hfill
  \includegraphics[height=0.3\textwidth]{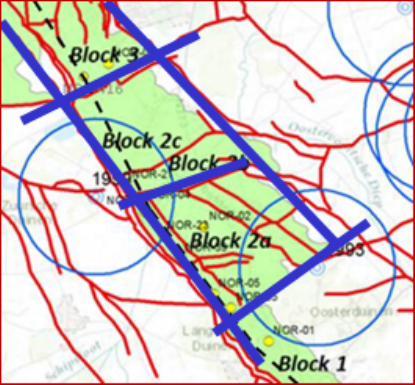}\hfill
  \hfill\null
  \caption{On the left: base Zechstein semblance map of the Norg UGS (in blue) and
surrounding area with traces of the bounding faults and localization of the recorded
seismic events \cite{NAM16}. On the right: conceptual map of the Norg field with major and minor faults
highlighted in blue and red, respectively.}
  \label{fig:geolMap}
\end{figure}

In order to identify the main mechanisms that govern fault reactivation in UGS fields, we use the modeling framework described above on
a simplified geological model representative of the typical features of
the Rotliegend UGS reservoirs, such as Norg and Grikpskerk \cite{NAM16}. These reservoirs
are bounded by normal faults with a significant throw (up to a 250~m) and consist of a few
compartments separated by internal faults (Fig.~\ref{fig:geolMap}). 
The gas fields are located between 2000 and 3000 m of depth, with the Rotliegend reservoir rock characterized by an
average net thickness of 150-200~m. Detailed information on the geological setting and
typical geometric features of UGS reservoirs in the Netherlands can be found, among
others, in these scientific papers \cite{Kra_etal13, Orl_etal13, fokker2013data,
wassing2017impact} and published reports \cite{TNO15, NAM16}.

\begin{figure}
  \centering
  \null\hfill
  \includegraphics[width=0.25\textwidth]{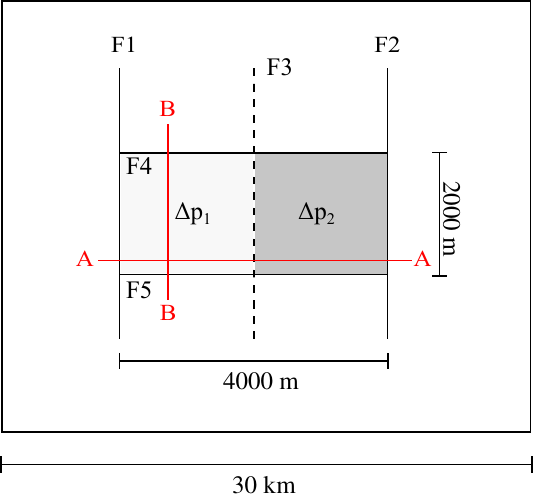}\hfill
  \includegraphics[width=0.70\textwidth]{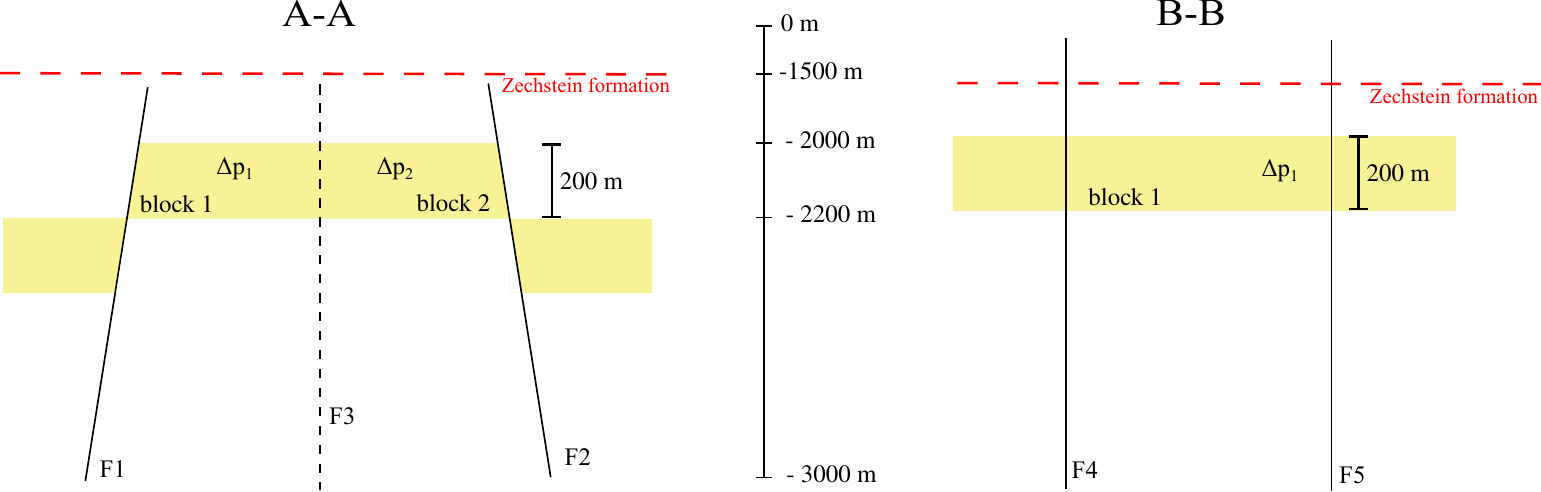}
  \hfill\null
  \caption{On the left: plain view of the model. On the right: vertical sections of the
conceptual model along the trace A-A and B-B shown on the left.}
  \label{fig:modelGeom}
\end{figure}

Based on those features, we define a conceptual model composed of two adjacent 
compartments, 2000$\times$2000~m wide, 200~m thick, and 2000~m deep, where UGS activities are carried out. 
The reservoir compartments are laterally confined by two
families of orthogonal faults, denoted as F1-F2 (parallel to $y$-axis) and F4-F5 (parallel to $x$-axis). Another fault, denoted as F3, separates the two reservoir blocks (Fig.~\ref{fig:modelGeom}).
The two compartments have only a partial hydraulic connection depending on the sealing properties of fault F3, so the pore pressure distribution in space and time may be different. Faults F1 and F2 are
inclined with respect to the vertical $z$-axis by a dip angle equal to $\pm 10^{\circ}$, while F3, F4 and F5 are vertical
faults, as shown in Fig.~\ref{fig:modelGeom}. The faults extend from -3000~m to -1600~m
depth, i.e., they terminate within the caprock sealing the reservoir, called Zeichestein
formation. Notice that the blocks have a 200-m offset along the vertical direction, corresponding to the entire
thickness of the reservoir, relative to the Rotliegend formation located in the sideburden. 

\subsection{Finite element-interface element discretization}
The reservoir is embedded in a 30-km wide square domain. The overall model size is much larger (about 10 times larger) than
the reservoir dimension to minimize the effects of the (arbitrary) boundary conditions on the solution
in the area of interest (Fig.~\ref{fig:modelGeom}) \cite{geertsma1973land}. The bottom of the model is 5000-m deep
and the land surface is located at the elevation of 0 m.

\begin{figure}
  \centering
  \null\hfill
  \includegraphics[width=0.45\textwidth]{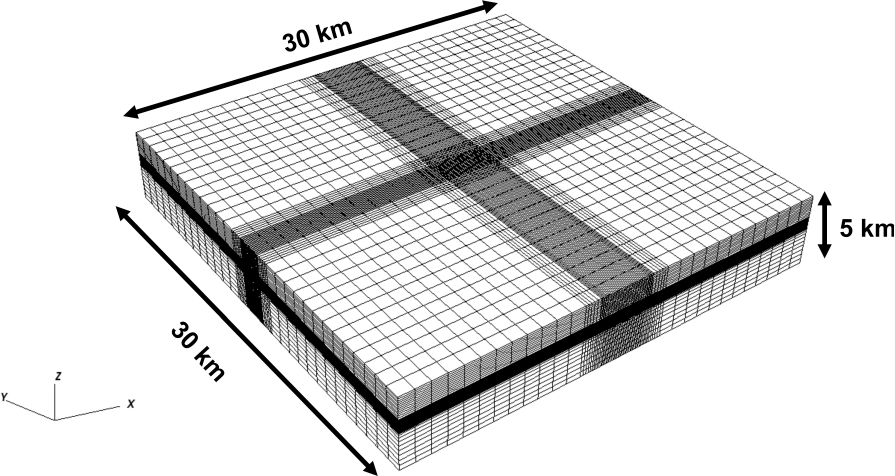}\hfill
  \includegraphics[width=0.45\textwidth]{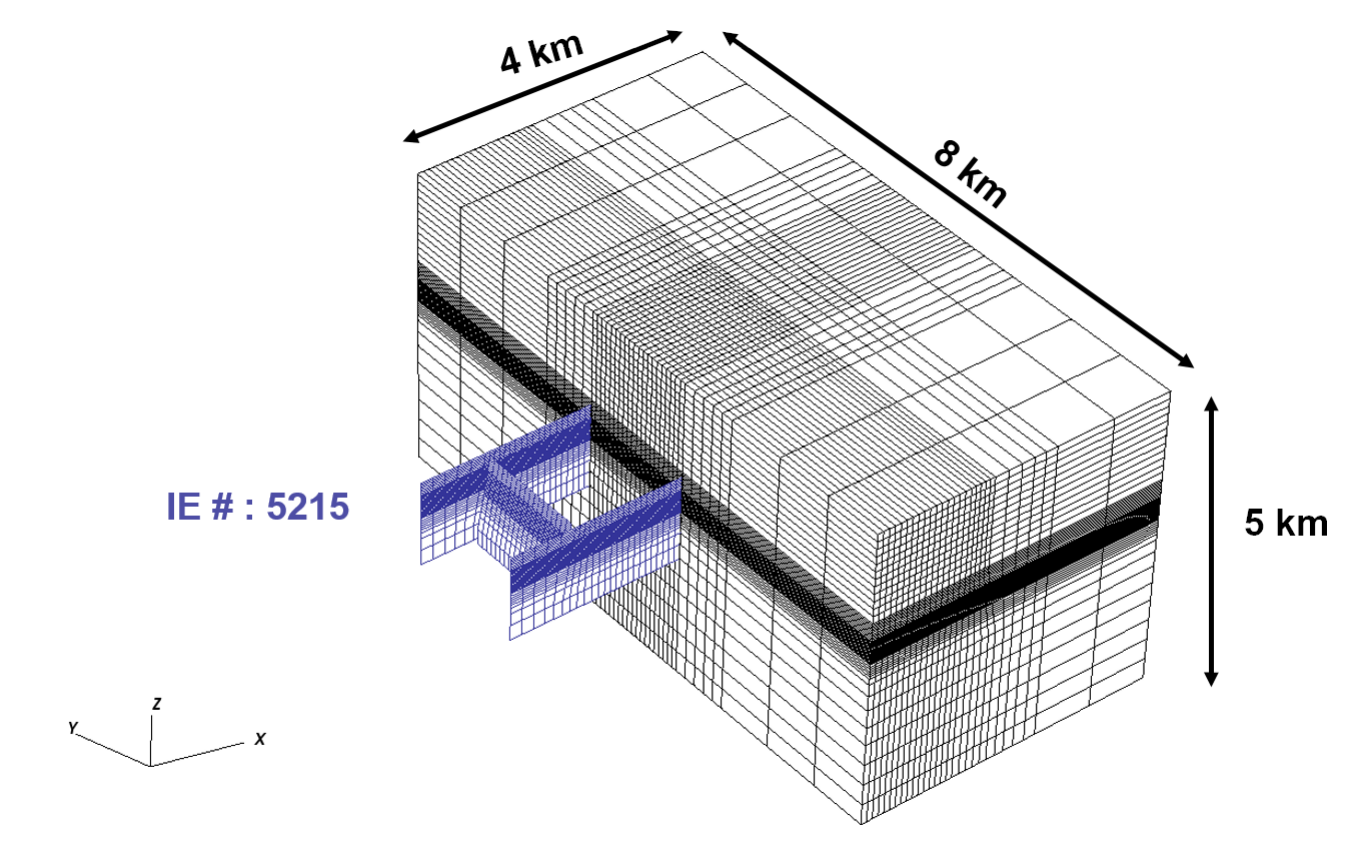}
  \hfill\null
  \caption{Axonometric view of the computational domain used for the geomechanical
simulations: full 3D finite element grid (left) and interface element grid (blue) embedded in a portion of the full 3D grid (right).}
  \label{fig:FEgrid}
\end{figure}

\begin{figure}
  \centering
  \includegraphics[width=0.50\textwidth]{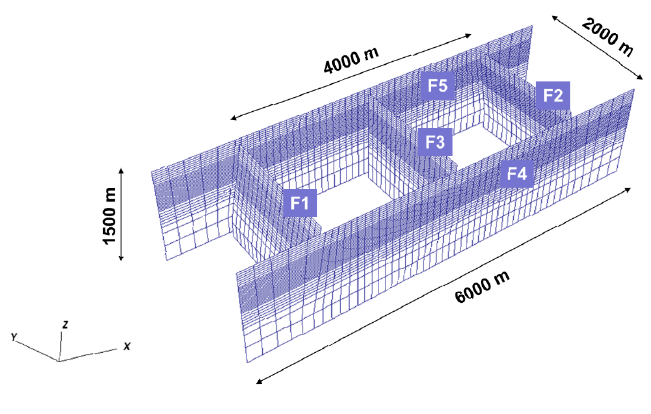} \hfill
   \includegraphics[width=0.40\textwidth]{figs/stressf0}
  \caption{On the left: interface element discretization of the fault discontinuities. The planar trace of F1, F2 and F3 is parallel to the y-axis, whereas that of F4 and F5 is parallel to the $x$-axis. F3 is the central fault separating the two reservoir compartments. On the right: initial normal stress with respect to the fault orientation. The principal stresses $\sigma_h$, $\sigma_H$ and $\sigma_v$ are parallel to the Cartesian axes. Faults F4 and F5 are more loaded because of their orthogonality to $\sigma_H$.}
  \label{fig:IEgrid}
\end{figure}

A 3D finite element mesh of the selected domain is built by
using hexahedral elements, which are particularly suitable for the selected symmetric
configuration with the faults parallel to the Cartesian axes. Fig.~\ref{fig:FEgrid} shows
an axonometric view of the full computational grid used in the geomechanical model. The mesh consists
of 253,165 nodes and 236,208 hexahedral elements with a finer discretization in the reservoir
layers, i.e., at depth between 2000 and 2200~m. The element size within the reservoir is
100$\times$100$\times$20~m. Fig.~\ref{fig:IEgrid} (left) shows the fault system embedded in the
continuous 3D grid as discretized by 5,215 interface elements. 
In the vertical direction, each fault is discretized by 40 elements, with an average size of 37.5~m.
The state of each element of the faults is synthetically evaluated with the aid of the \textit{criticality index} $\chi$ defined as:
\begin{equation}
  \chi = \frac{\|\vec{t}_T\|_2}{\tau_{\max}} = \frac{\|\vec{t}_T\|_2}{c - t_N\tan
\left( \varphi\left(\|\vec{g}_T\|_2\right)\right)}.
  \label{eq:chi}
\end{equation}
From Eq.~\eqref{eq:chi}, it is easy to see that $\chi \in [0,1]$, where 0 is associated
with the safest condition and 1 to plastic sliding.
We would like to clarify that we intentionally avoided utilizing domain symmetry and
modeling for only half of the domain to be able to simulate asymmetric loading
conditions as well.
\af{Finally, we emphasize that in our model the different faults end close to each other, but do not intersect. Fault intersections are deliberately not modeled at this development stage of our code, as their inclusion would introduce an important implementation burden without contributing in a significant way to the overall physical description of the process.}

\begin{figure}
  \centering
  \includegraphics[width=0.40\textwidth]{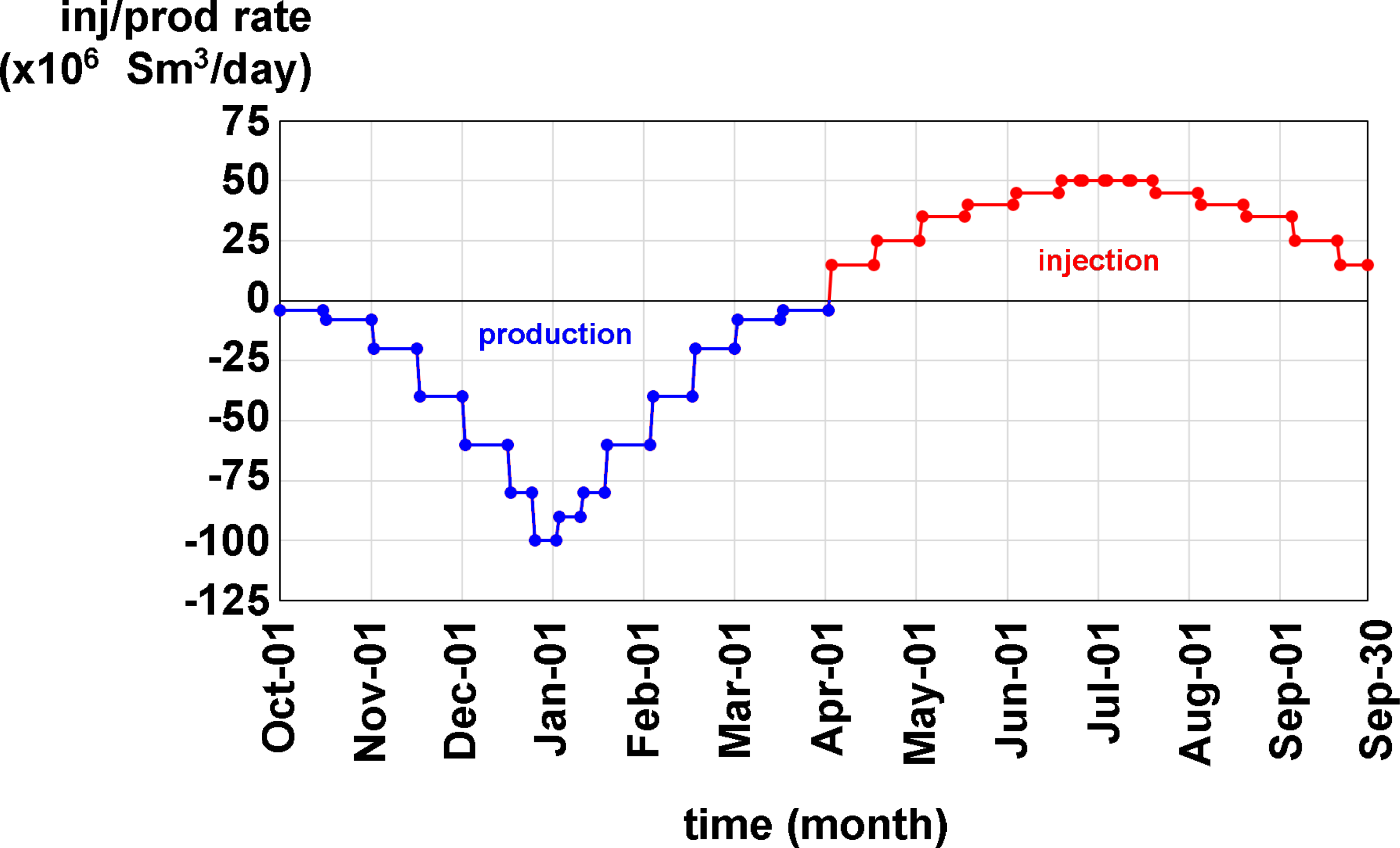}
  \caption{Time behavior of the production/injection rate used in the OPM simulation.}
  \label{fig:OPMrate}
\end{figure}

\subsection{Pore pressure variation}
\begin{figure}
  \centering
  \includegraphics[width=0.70\textwidth]{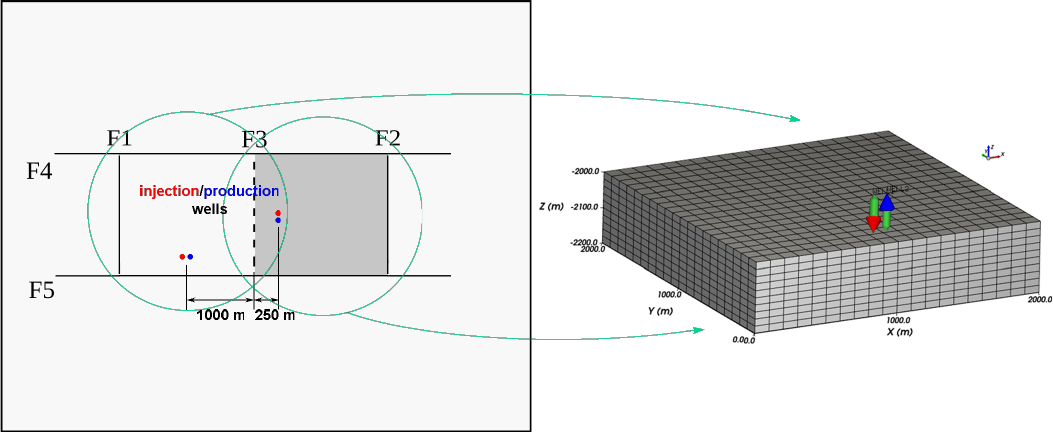}
  \caption{Location of the injection/production wells in the two reservoir
blocks (left) and axonometric view of the 3D computational grid used in OPM to simulate the
injection/production phase in each reservoir compartment (right). The OPM mesh exactly corresponds to the finite element grid of a single
block within the 3D geomechanical model.}
  \label{fig:OPMgrid}
\end{figure}

As previously mentioned, in this work we adopt a one-way coupled approach, thus the
multiphase flow prediction is computed first. The simulation is
performed through the open-source reservoir simulator OPM Flow \cite{OPM,
rasmussen2021open}. As a reference scenario, a typical year-long cycle of UGS activity has
been considered, with the injection-production history represented in
Fig.~\ref{fig:OPMrate}. Fig.~\ref{fig:OPMgrid} shows the location of the
injection/production wells with respect to the fault system. Note that, to avoid any
interpolation among computational grids, the OPM finite volume mesh exactly corresponds to
the finite element grid of a single block within the 3D geomechanical model.
The characteristic horizontal and vertical permeability of a reservoir in the study area
are $k_h = 600$~mD and $k_v = 300$~mD, respectively. They are kept constant during the
simulation. The ``working gas'' volume amounts to 6.5$\times$10$^9$~Sm$^3$ per
compartment.
\af{In these simulations, we do not account for capillary pressure in order to simplify the modeling process.}

\begin{figure}
  \centering
  \includegraphics[width=0.40\textwidth]{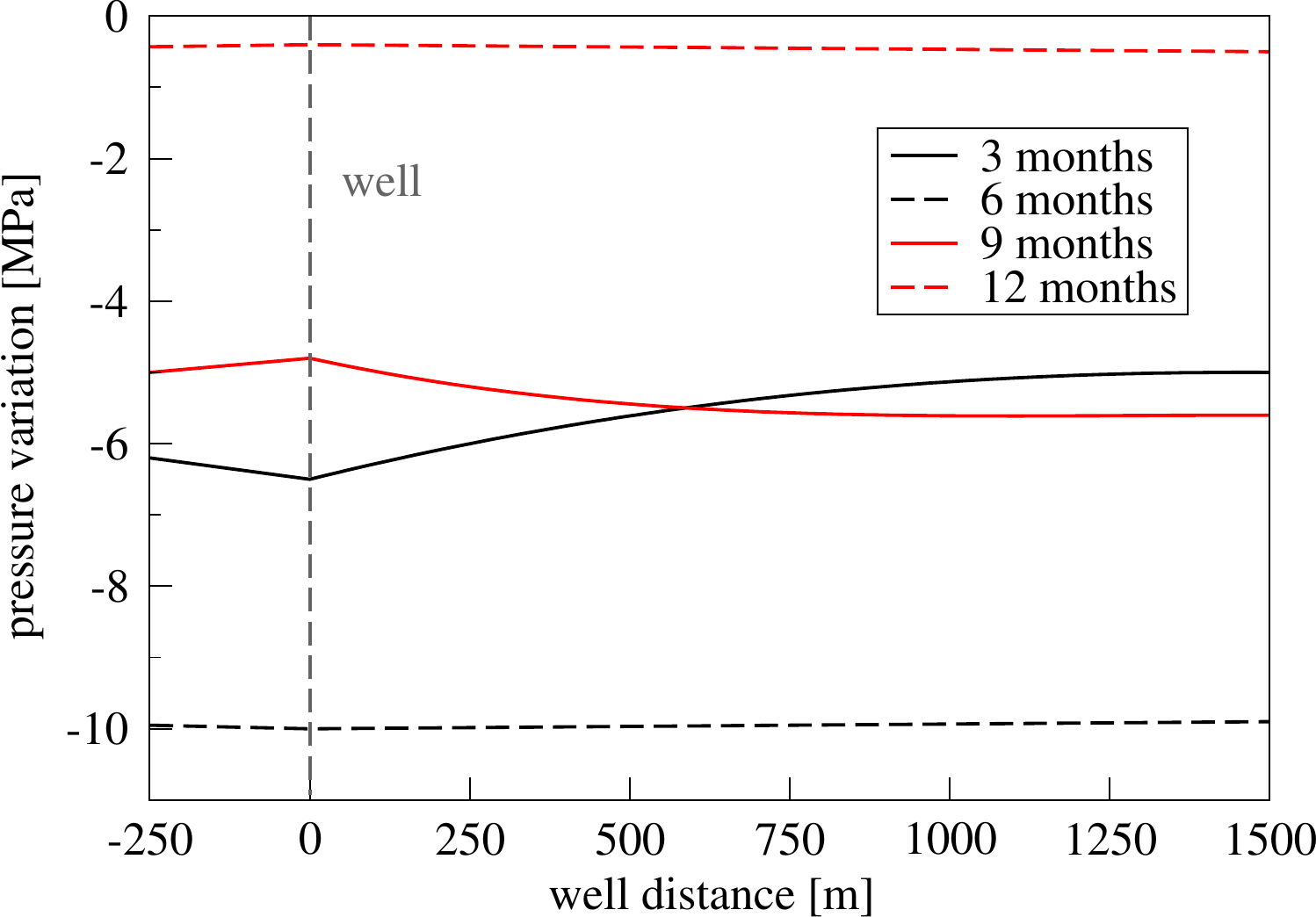} \hspace{1truecm}
  \includegraphics[width=0.40\textwidth]{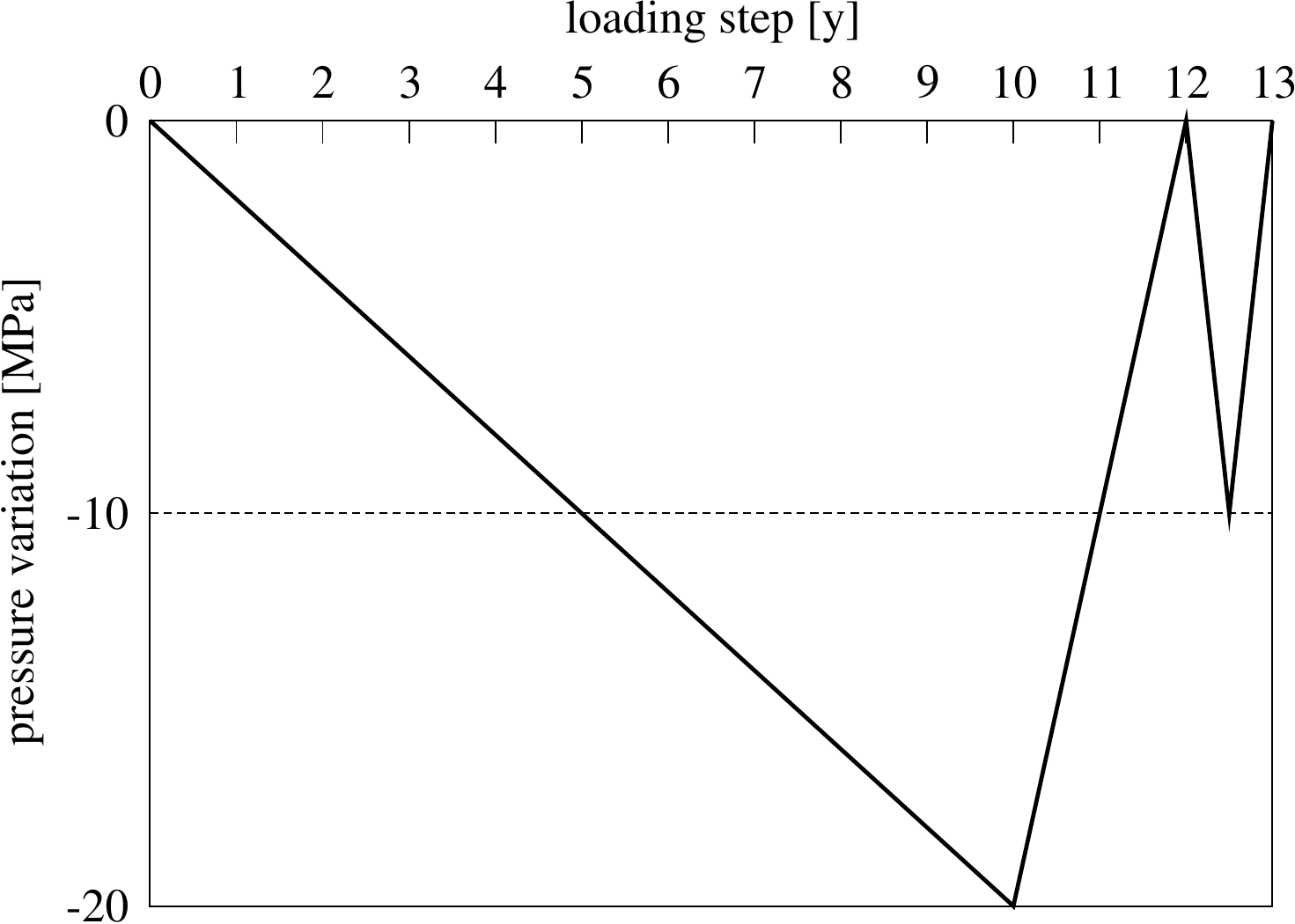}
  \caption{On the left: depth-averaged value of the pore pressure variation during a production/injection cycle as obtained by the OPM flow simulator. On the right: sketch of the pore pressure variation over time as prescribed within the reservoir compartments.}
  \label{fig:OPMpress}
\end{figure}

The OPM Flow results in terms of pressure variation are summarized in Fig.~\ref{fig:OPMpress} (left). The figure shows
the depth-averaged pressure behavior along a section passing through the
production/injection wells every 3 months. After 3 months the maximum production rate is
achieved, after 6 months the production phase ends, after 9 months the maximum injection
rate is met, and finally, after 12 months the simulation ends. Notice that the pressure
perturbation during the entire production (or injection) phase is almost uniform in space and varies approximately within the interval between 0 and -10~MPa with respect to the initial value $p_i$.
This outcome shows that, in agreement with previous modeling studies \cite{Bui_etal19}, the spatial gradient of the pore pressure variation into each compartment of Rotliegend reservoirs is expected to be quite limited. Hence, considering a constant pressure variation value for each reservoir block appears to be a reasonable assumption. 

Based on a typical production file of Dutch UGS reservoirs \cite{MIT09, NAM16}, the entire pressure history prescribed in each compartment is sketched in Fig.~\ref{fig:OPMpress} (right). We assume a 10-y duration for the PP phase, where the pressure drops linearly by up to
20~MPa. After this period, a 2-year CGI
phase follows, where the pressure recovers to the initial (undisturbed) value
$p_i$, and then the UGS cycles start. They are characterized by a 6-month extraction
period, during which the pressure drops by 10~MPa, and a 6-month
injection period, when the pressure returns to $p_i$ in agreement with the outcome of OPM flow.

\subsection{Simulated scenarios}
To evaluate the capabilities of the presented numerical model and understand the possible
mechanisms causing fault reactivation during CGI and UGS, a few scenarios are simulated in
the typical setting of the Rotliegend reservoirs in the Netherlands. The main geological and geomechanical parameters are reported in Tab.~\ref{tab:params}. For the sake of simplicity, a linear elastic behavior is assumed in the reservoir during the UGS activities.

Standard conditions with zero
displacement and zero pore pressure variation on the outer and bottom boundaries are prescribed, whereas the land surface is a traction-free boundary.


\begin{table}
\caption{Formation-dependent geomechanical parameters. See Fig.~\ref{fig:modelGeom} for a
detail on the depths.}
\label{tab:params}
\centering
\begin{tabular}{l|c|c|c} \hline
  layer & density [kg/m$^3$] & Young modulus [GPa] & Poisson ratio \\
  \hline
  Overburden & 2200 & 10.0 & 0.25 \\
  Upper Zechstein Salt (-1500 to -1800~m) & 2100 & 35.0 & 0.30 \\
  Lower Zechstein Salt (below -1800~m) & 2100 & 20.0 & 0.30 \\
  Reservoir (Upper Rotliegend) & 2400 & 11.0 & 0.15 \\
  Underburden & 2600 & 30.0 & 0.20 \\ \hline
\end{tabular}
\end{table}

To initialize the simulation, the undisturbed stress regime must be prescribed. We assume the initial pressure regime following a hydrostatic distribution and the principal effective
stress tensor directions aligned with the Cartesian axes. In particular, $\sigma_1 =
\sigma_v = \sigma_z$, $\sigma_2 = \sigma_H = \sigma_y$, and $\sigma_3 = \sigma_h =
\sigma_x$, where $\sigma_v$ denotes the vertical compressive stress, and $\sigma_H$ and $\sigma_h$ the largest and smallest compressive horizontal principal stress, respectively. At the reservoir average depth, i.e., $z=-2100$~m, we have $\sigma_v = -25.4$~MPa,
$\sigma_h = M_1 \sigma_v = -18.8$~MPa, $\sigma_H = M_2 \sigma_v = -21.1$~MPa, with $M_1 =
0.40$ and $M_2 = 0.47$. The initial normal stress acting on the faults is shown in
Fig.~\ref{fig:IEgrid} (right). We emphasize that the initial stresses on the faults are only due to gravity contributions.

Two scenarios have been simulated based on the parameters describing the Coulomb
frictional criterion. In the reference scenario (scenario~1), $\varphi_s = 30^{\circ}$ and
fault weakening is not accounted for. The effect of slip-weakening behavior is investigated in
scenario~2, where the friction angle reduces from $\varphi_s = 30^{\circ}$ to $\varphi_d =
10^{\circ}$ in a slip distance of $D_c = 2$~mm. Cohesion $c = 2$~MPa in both scenarios.
We acknowledge that considering the same properties for each fault could be an
oversimplification, but for the purpose of analyzing an idealized case, we accept this
assumption.
Finally, the time step is 1~year during the PP phase, then during the CGI and UGS phases it is reduced to 15~days (a 2~months time step is used to show the results).

\subsection{Numerical performance of the computational contact model}
\begin{figure}
  \centering
  \includegraphics[width=0.40\textwidth]{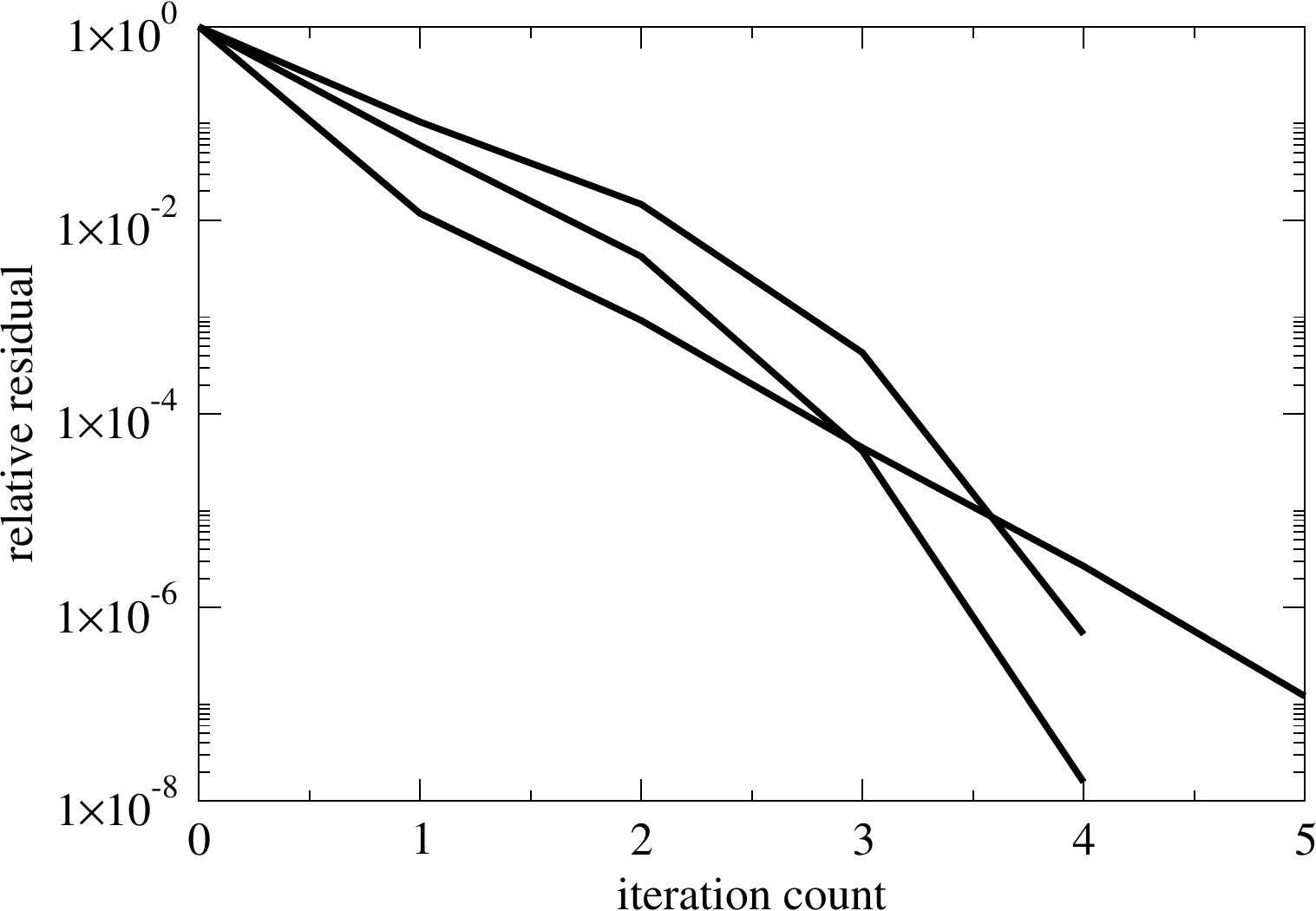}\qquad
  \includegraphics[width=0.41\textwidth]{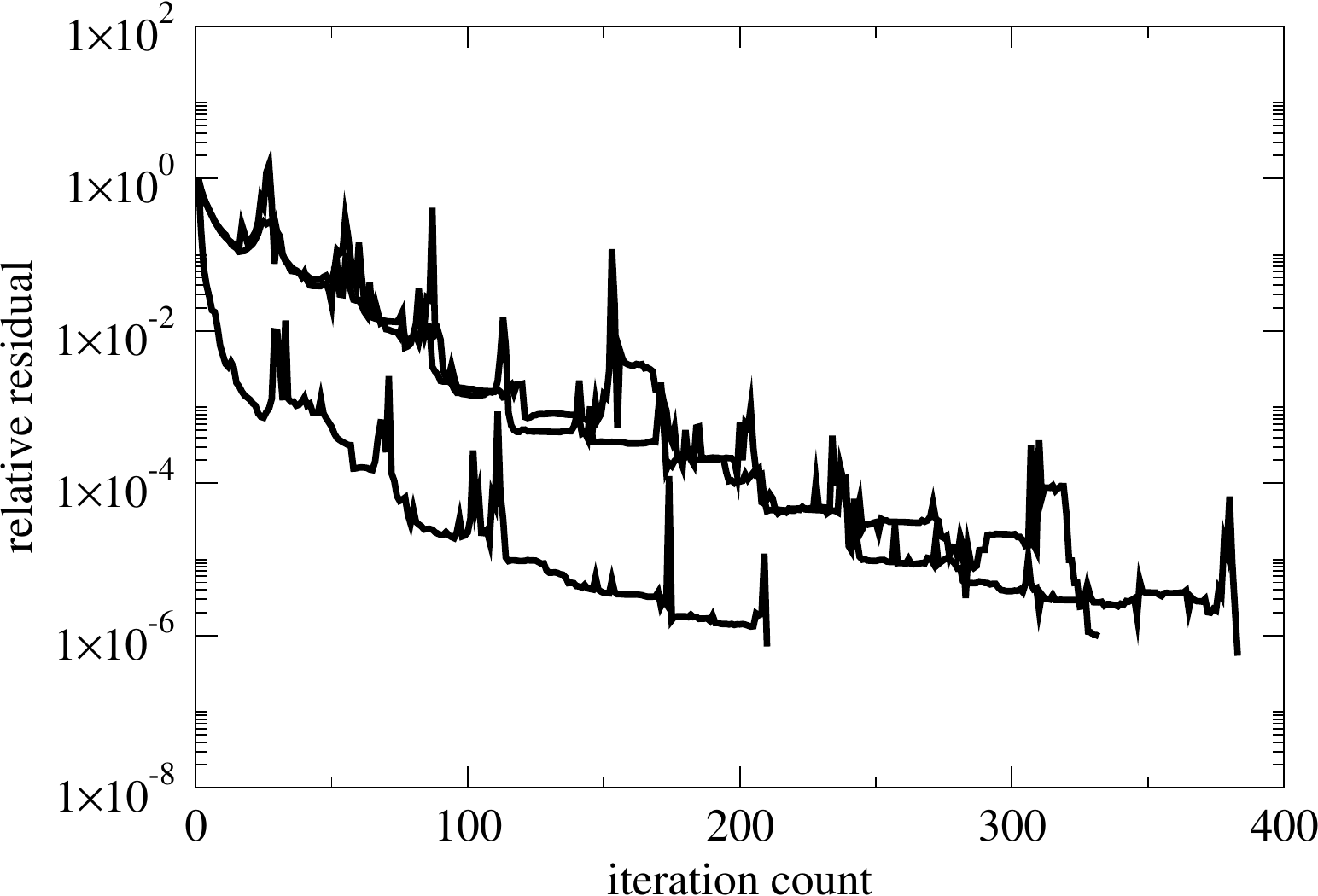}
  \caption{\af{Non-linear convergence profile for some representative active-set iterations with fault activation occurrences (left). Convergence profile of the preconditioned Bi-CGStab solver for some representative inner linear system with the Jacobian matrix.}}
  \label{fig:convProf}
\end{figure}

\af{
At each time-step, 
the active-set strategy identifies the active and inactive Lagrange multipliers, then a Newton iteration solving the constrained non-linear problem is applied. If the consistency check at convergence of the non-linear iteration is satisfied, the simulation moves at the following time-step, otherwise the new set of active and inactive multipliers is identified.
%
If there is no fault reactivation ($\chi < 1$ everywhere), all Lagrange multipliers are unknown and the overall problem is simply linear. 
This situation occurs on most of the time steps, as we will see in detail in the next subsection.
In this case, the active-set strategy and the Newton method converge in a single iteration. 
When the condition $\chi = 1$ is attained in at least one of the elements discretizing the faults, the situation changes. Given the nature of the application, i.e., geological faults at the macro-scale where sliding occurs in a relatively small area with respect to the size of the entire fault surface, the number of active-set steps is limited, totaling no more than 4 in our simulations. Similarly, the number of Newton iterations required to solve the nonlinearity introduced by the friction law is also small, not larger than 7. A few representative convergence profiles of the Newton method are shown in Fig. \ref{fig:convProf} (left), which provide evidence of the expected quadratic order in asymptotic conditions.} 

\af{
At each iteration of the Newton method, a non-symmetric generalized saddle-point linear system with global size equal to 775,140 has to be solved.
For this task, we employ
an ad-hoc preconditioned Bi-CGStab algorithm. The use of preconditioners discussed in the
literature \cite{franceschini2019block, franceschini2022reverse},
specifically the Reverse Augmented Constrained Preconditioner (RACP) approach,
allows for a robust convergence in all system solutions with a relatively low number of iterations, typically less than 400. A few representative Bi-CGStab convergence profiles are provided in Fig. \ref{fig:convProf} (right).} 

\af{On summary, advancing in time requires 2 linear system solutions with no fault activation (linear problem), or a variable number of system solutions, generally between 5 and 20, if sliding occurs. Each linear system solution, on a standard workstation equipped with a CPU at 2.4 GHz and 64 Gbyte RAM, requires on average 80 seconds to set-up the preconditioner and 140 seconds to iterate until convergence, resulting in a total time of approximately 220 seconds per linear system. This computational efficiency allows to run several simulations in an acceptable amount of time with no specific hardware requirements.}

\subsection{Hazard of fault reactivation}
The objective of the representative simulations reported herein is to evaluate the fault reactivation hazard during the different stages of the UGS activities in the conceptual reservoir. 
For this reason,
we mainly focus on the criticality index $\chi$ defined in Eq.~\eqref{eq:chi}. For the sake of clarity and ease of readability, $\chi$ is represented for each fault as a
function of depth only, i.e., for each $z$-value we compute the $\chi$ average for the stripe of interface elements located at the same depth.
Another significant quantity is the
maximum sliding, i.e., the maximum value of $\|\vec{g}_T\|_2$ simulated along each fault. These two
quantities are closely related each other, since a single element can slide only when $\chi = 1$.
However, we prefer to propose an averaged version of $\chi$, so as to obtain information on the criticality
state of the entire fracture at a given depth.

The last quantity used to interpret the
results and analyze the fault behavior is the tangential component of the traction. In particular, we use $t_{T,z}$, i.e., the vertical
component of $\vec{t}_T$. Usually, the 2-norm of the tangential traction is
analyzed, i.e., $\|\vec{t}_T\|_2$, but this does not provide information on the shear direction. However, thanks to the symmetric geometry of the conceptual model, in some locations there
is no horizontal component of $\vec{t}_T$, thus, $\|\vec{t}_T\|_2 = |t_{T,z}|$. The two
quantities share the same modulus, but the vertical component carries additional
information on the sliding direction.

\begin{figure}
  \centering
  \includegraphics[width=\textwidth]{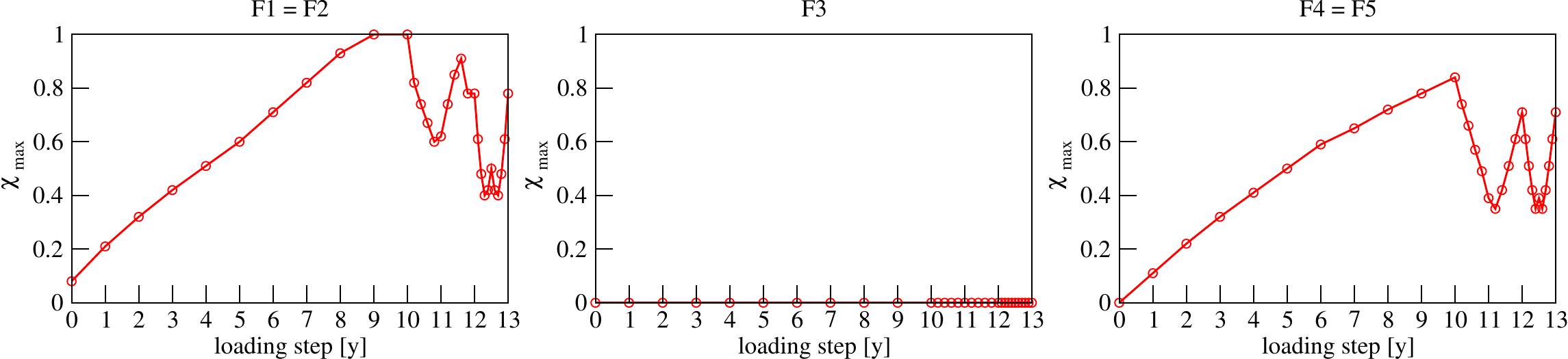}
  \caption{Behavior of $\chi_{\max}$ from all the loading steps for each fault. Note that
due to symmetry F1 and F2 behave identically, as well as F4 and F5.}
  \label{fig:chiAll}
\end{figure}

\begin{figure}
  \centering
  \null\hfill
  \includegraphics[align=c,width=0.30\textwidth]{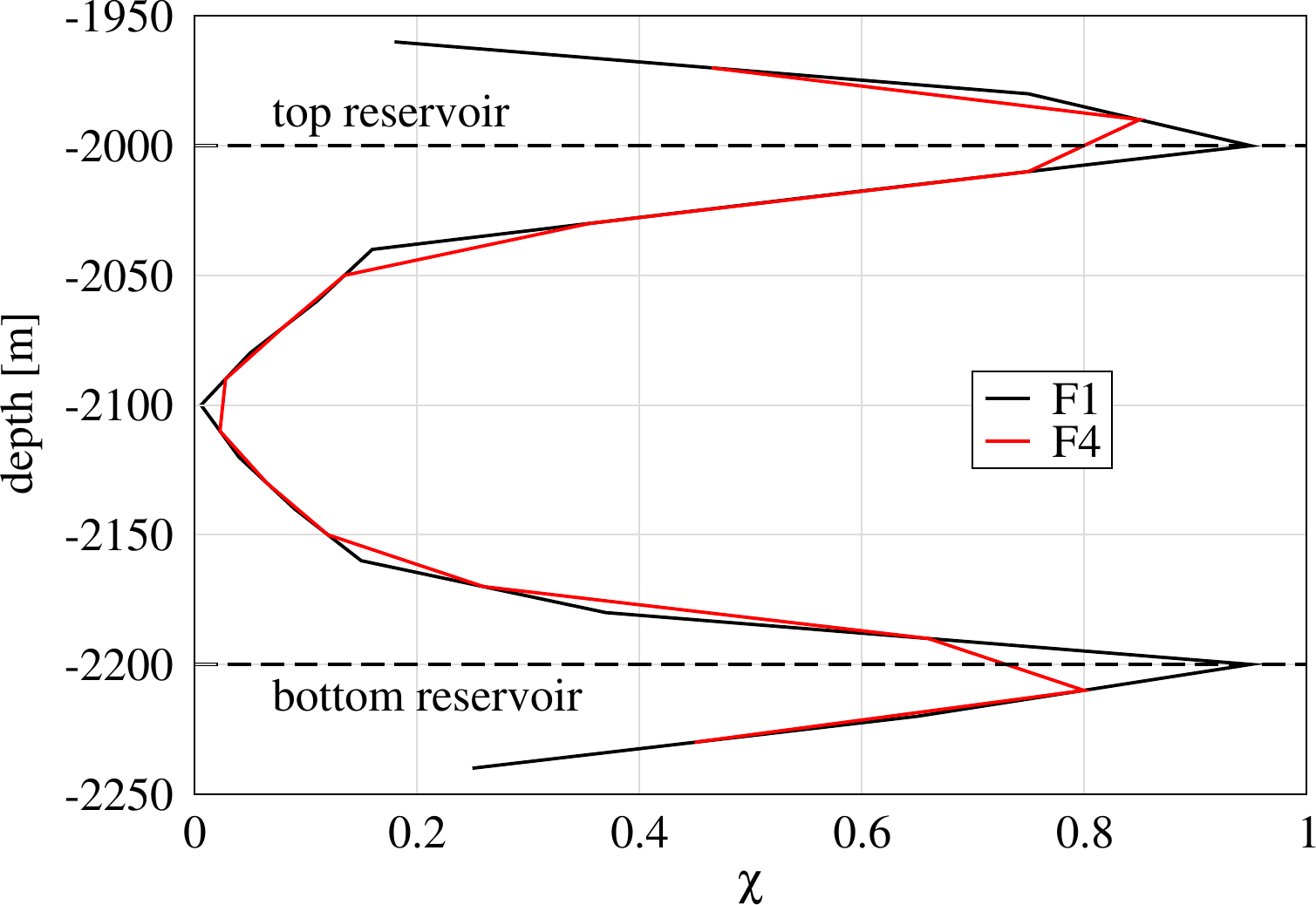}\hfill
  \includegraphics[align=c,width=0.40\textwidth]{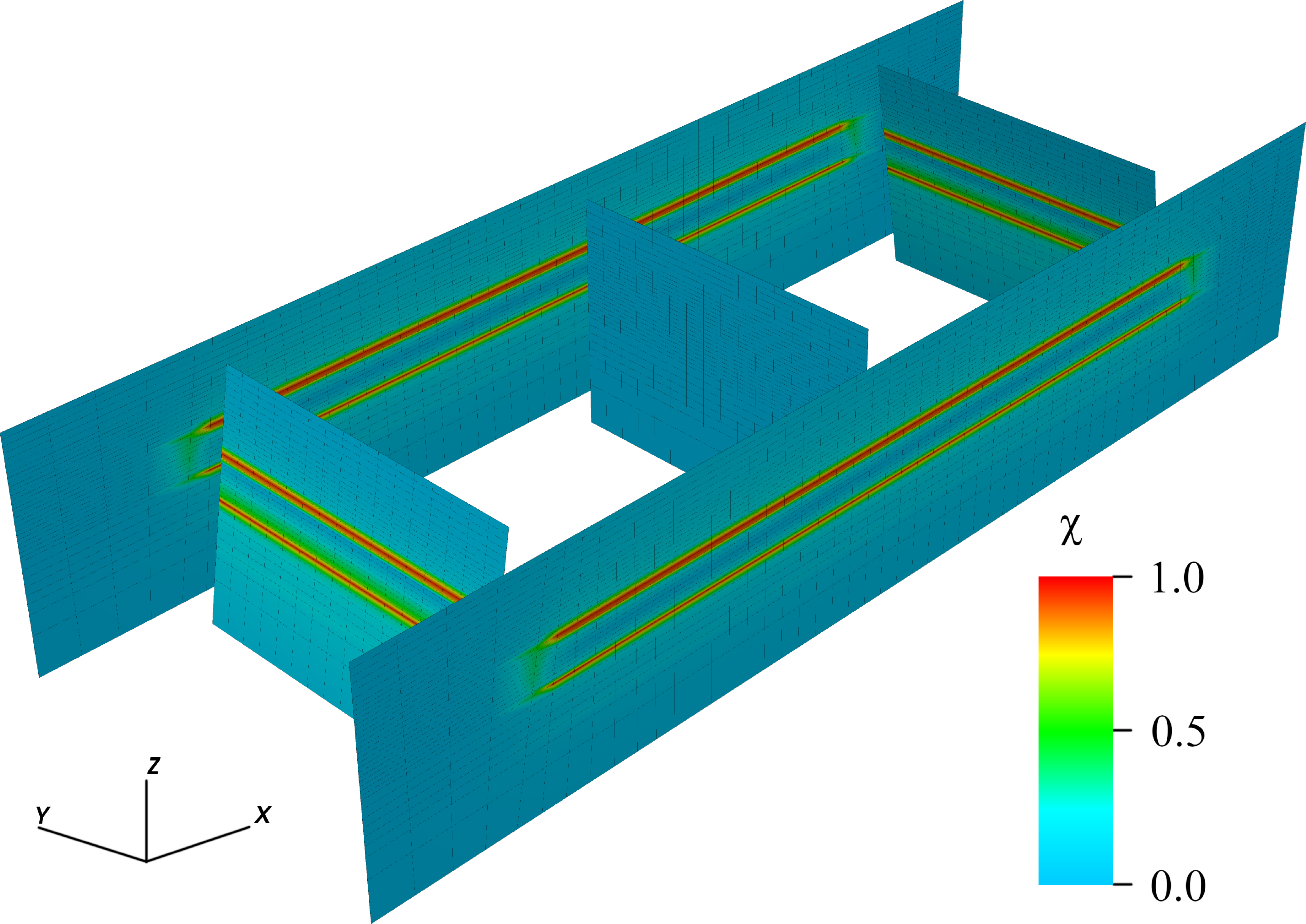}
  \hfill\null
  \caption{On the left: behavior of the criticality index $\chi$ vs depth at loading step
10 (end of PP) for faults F1 and F4. On the right: $\chi$ factor on all the fracture
surfaces at loading step 10.}
  \label{fig:chiz}
\end{figure}

The mechanisms for the possible fault reactivation have been investigated in scenario~1. The value
$\chi_{\max} = 1$ is reached on faults F1 and F2 at loading step 9, with $\chi_{\max}$ up
to 0.8 at the end of CG and UGS injection phases (Fig.~\ref{fig:chiAll}). Conversely,
$\chi = 0$ on fault F3 irrespective of the loading step due to the symmetry of the
geometry and loading configurations. A comparison between the behavior versus depth of the
criticality index along fault F1 and F4 and the distribution of $\chi$ on the whole fault
system at the end of PP are shown in Fig.~\ref{fig:chiz}. Notice that the most critical
condition develops along the top and bottom of the reservoir in agreement with previous
modeling study \cite{Hau_etal18}.
Moreover, faults F4 and F5 exhibit smaller values of $\chi$ with respect to F1 and F2, showing that a sub-vertical orientation is usually more likely to reactivate.

\begin{figure}
  \centering
  \null\hfill
  \includegraphics[align=c,width=0.3\textwidth]{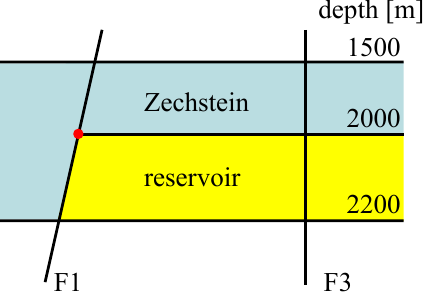}\hfill
  \includegraphics[align=c,width=0.38\textwidth]{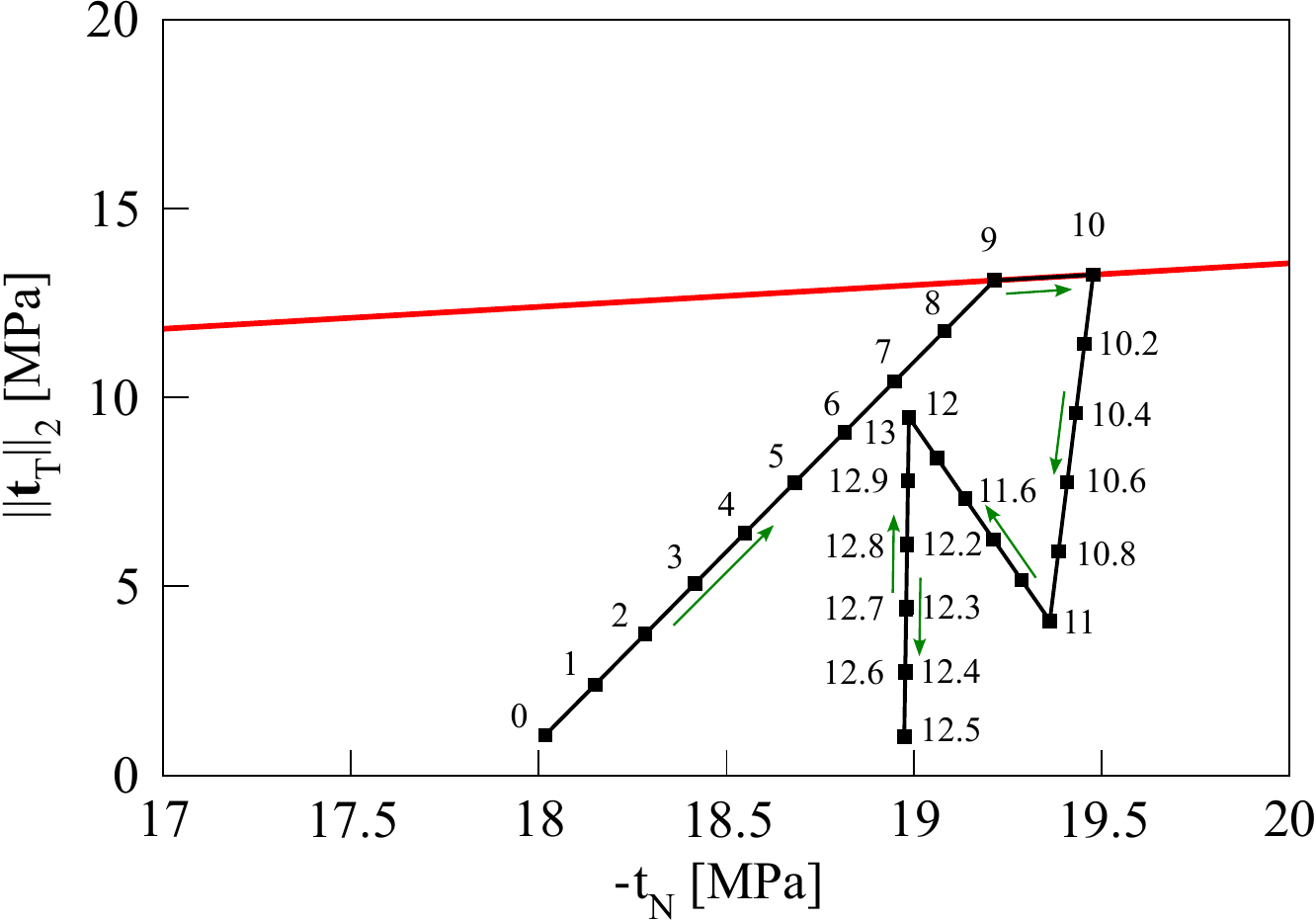}
  \hfill\null
  \caption{On the left: location of the selected element. On the right: stress path
$\|\vec{t}_T\|_2$ vs $-t_N$ for the element highlighted on the left sketch by a red dot. The red
line is the yield bound. Numbers along the path denote the loading steps. It can be easily
recognized the primary production (loading steps 1 to 10), the cushion gas injection (loading steps
10 to 12) and the underground gas storage (loading steps 12 to 12.5 -- production -- and 12.5 to 13 -- injection).}
  \label{fig:stressPath}
\end{figure}

Fig.~\ref{fig:stressPath} shows the stress path in the $t_N-\|\vec{t}_T\|_2$ plane experienced by a representative element located on
fault F1 at the top of the reservoir. The actual stress state touches the yield bound at the
loading step number 9 and remains on the yield surface until the end of the PP (loading step 10). During CGI, the stress state initially departs from the
yield condition but returns close to it during the last part of the injection when the
pressure recovers to the initial value. UGS behaves elastically on a new path compared to what was experienced during the last part of the CGI phase, with an almost constant $t_N$ value. Again, the stress state approaches the critical condition at the end of the UGS injection phase when $p$ returns to $p_i$.

\begin{figure}
  \centering
  \null\hfill
  \includegraphics[align=c,width=0.35\textwidth]{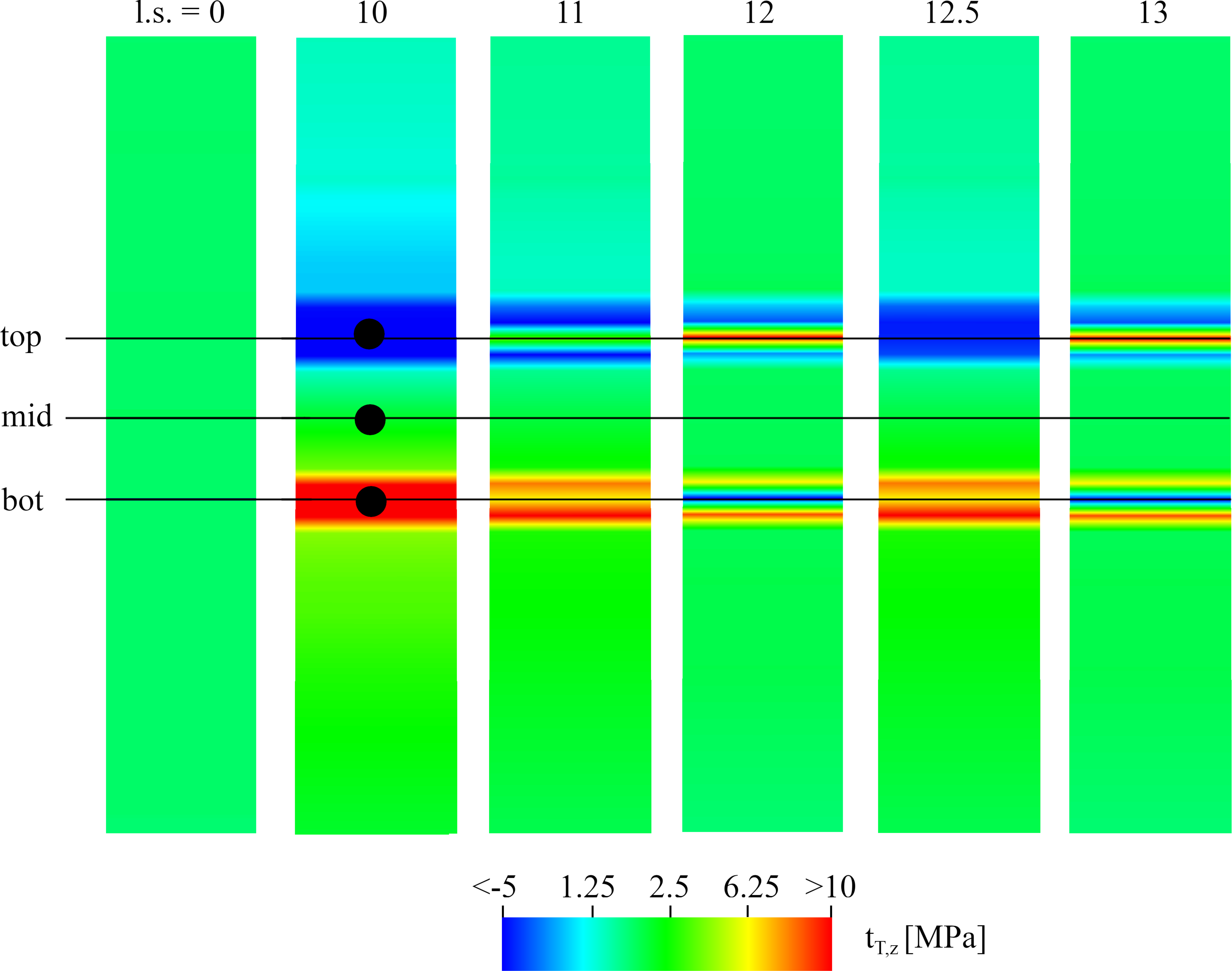}\hfill
  \includegraphics[align=c,width=0.35\textwidth]{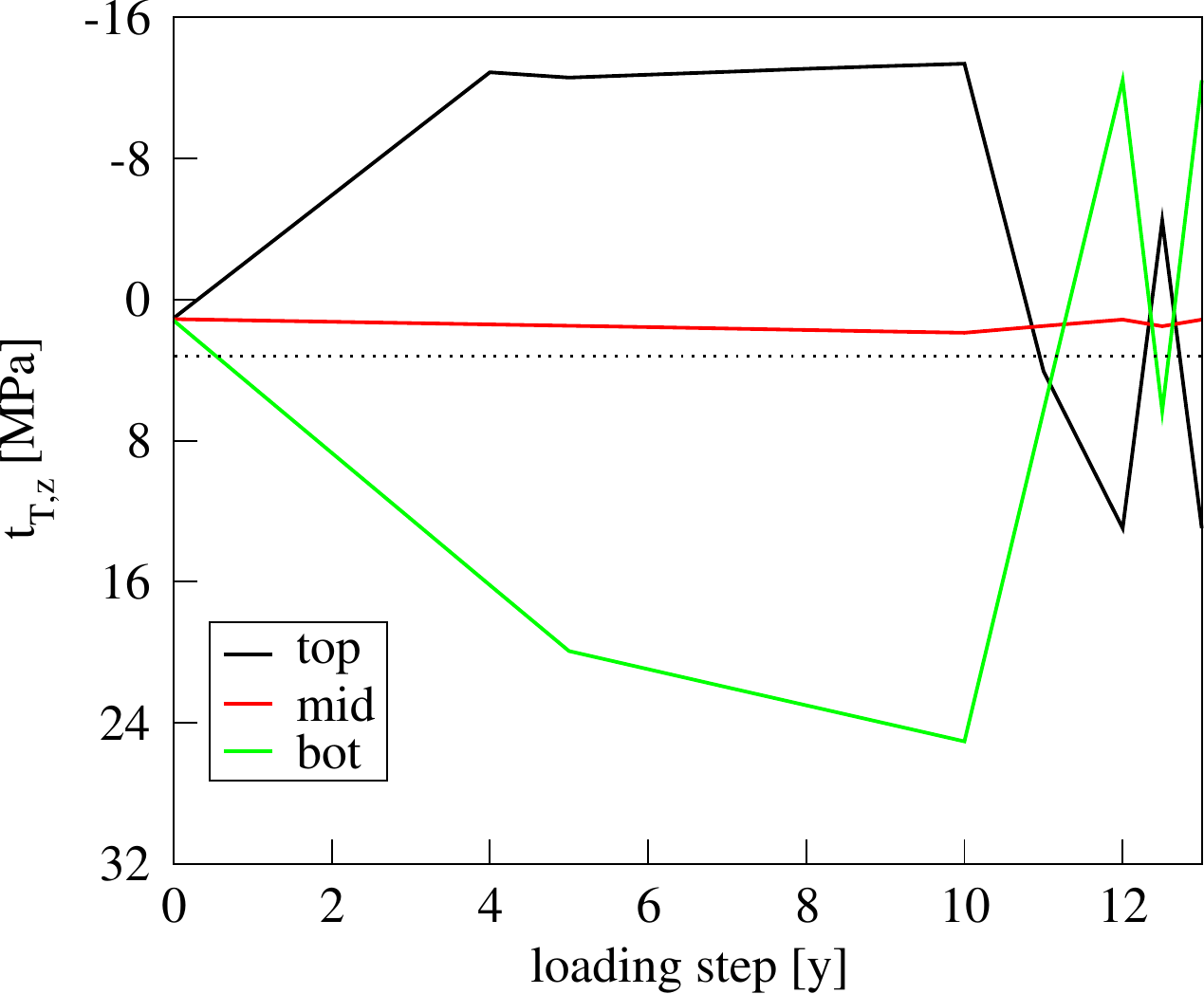}
  \hfill\null
  \vspace{10pt}
  \includegraphics[align=c,width=\textwidth]{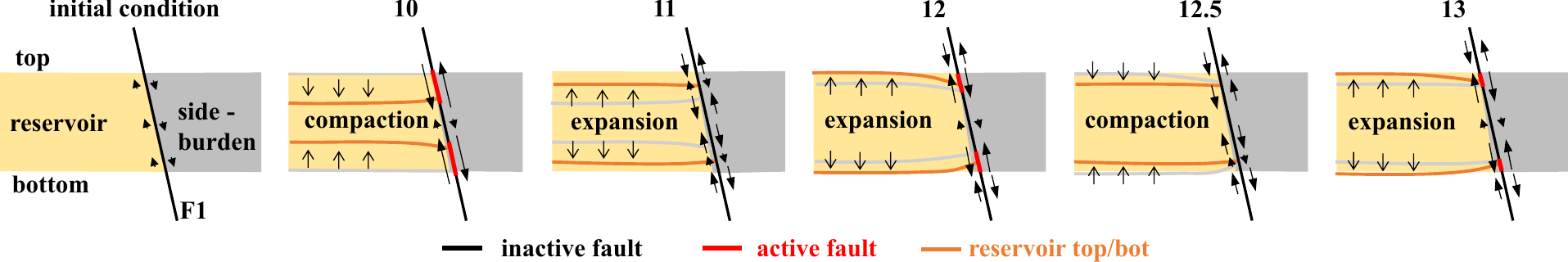}
  \caption{On the top left: distribution of vertical component of the shear stress $t_{T,z}$
for the loading steps (l.s.) 0 (initial condition), 10, 11, 12, 12.5, 13 on fault F1 (dip
= 10$^{\circ}$). On the top right: time behavior of $\|\vec{t}_T\|_2$ for the points denoted by the thick black dots in the previous frame located at the top, bottom, and center of the reservoir. Positive values mean that the shear
stress is directed upward. On the bottom: sketches
representing the reservoir deformation, shear stress direction, and inactive/active
portions of the fault at the same loading steps.}
  \label{fig:tau}
\end{figure}

A deeper explanation for this behavior can be found by analyzing the actual direction of the shear stress. Fig.~\ref{fig:tau} shows the vertical component of tangential traction $t_{T,z}$ on fault F1 at loading steps 0 (initial condition), 10, 11, 12, 12.5, and 13.
This component is meaningful because of the symmetry of the model, indeed, we have that
$\|\vec{t}_T\|_2 = |t_{T,z}|$. Sketches of the reservoir-fault-sideburden conditions are
provided for the same loading steps. The initial shear stress differs from the null value
because of the fault dip. The highest value of $t_{T,z}$ is observed at the end of PP (loading
step 10). The reservoir compaction induced by pressure depletion is accompanied by fault
reactivation. Note that positive and negative shear stress characterize the bottom and top of the reservoir, respectively. As physically expected because of the compaction mechanism, the direction of the shear stress is oriented toward the
center of the reservoir. When CGI starts, the shear stress orientation changes and the reactivated part of the fault returns stick.
At loading step 11, half of the pore pressure change has been recovered. As the reservoir
expands due to pressure recovery, $t_{T,z}$ decreases on the top and bottom of the reservoir (the
orientation remains the same but the absolute value decreases) and an almost null
$t_{T,z}$ is obtained at this step in the previously sliding IEs. Unlikely, $t_{T,z}$
does not change significantly for the elements surrounding the activated stripes of the fault. The reservoir
continues to recover pressure and re-expand until loading step 12. During this second part
of CGI shear stress increases, with a sign opposite to that experienced during PP
(Fig.~\ref{fig:tau}). A mirror behavior occurs for the IEs at the reservoir bottom.
Therefore, expansion during CGI increases the criticality condition of the fault (mainly
at the reservoir top and bottom) due the stress re-distribution after the sliding
developed over the PP. Fig.~\ref{fig:chiAll} shows that faults F1 and F2 approach the
criticality state ($\chi_{\max} > 0.8$) when the pressure recovers the initial value at the
end of CGI and UGS injection phase, i.e., in a pressure state close to the initial undisturbed one, which is not generally expected to be associated to fault reactivation.

\begin{figure}
  \centering
  \null\hfill
  \includegraphics[width=0.38\textwidth]{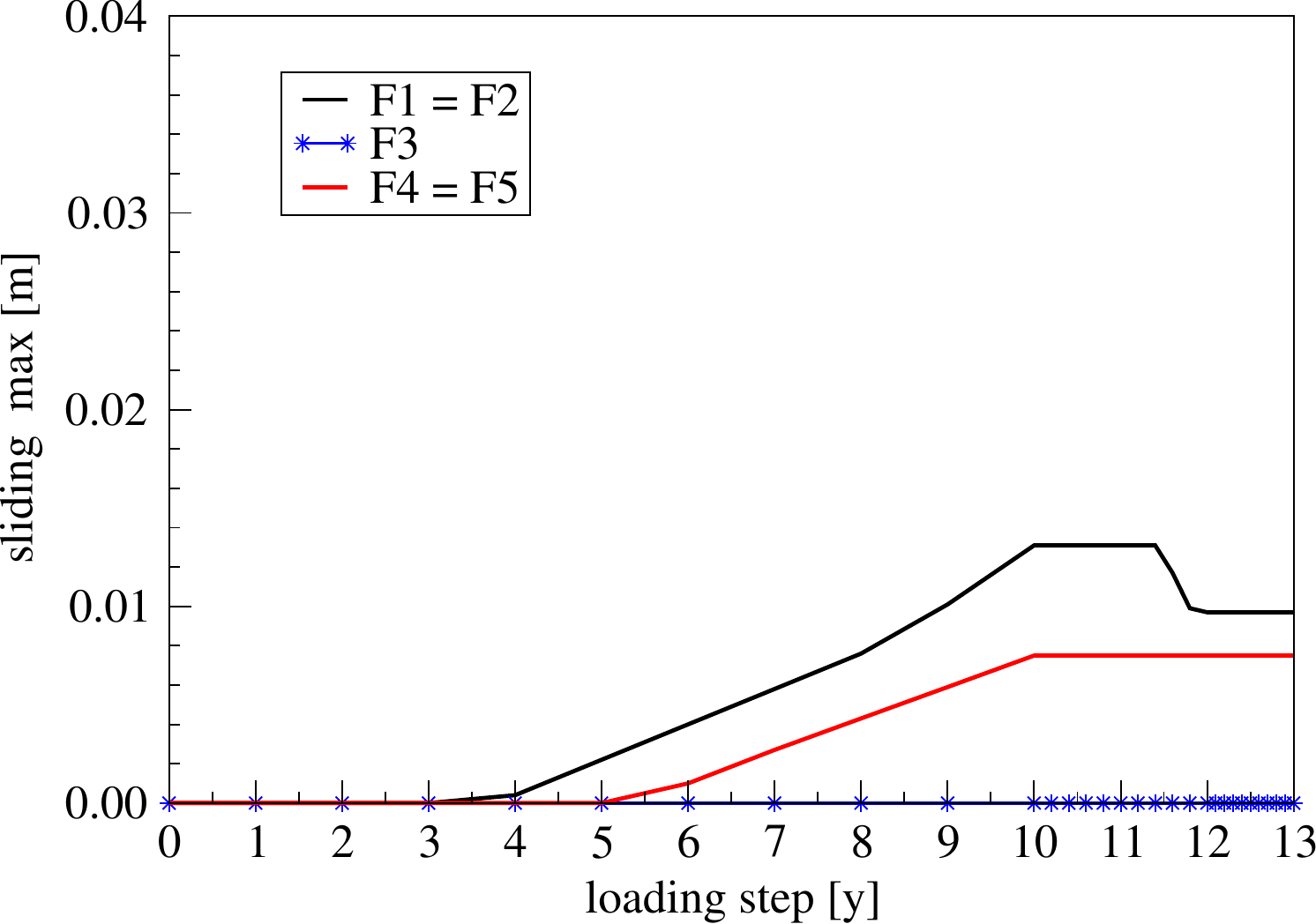}\hfill
  \includegraphics[width=0.38\textwidth]{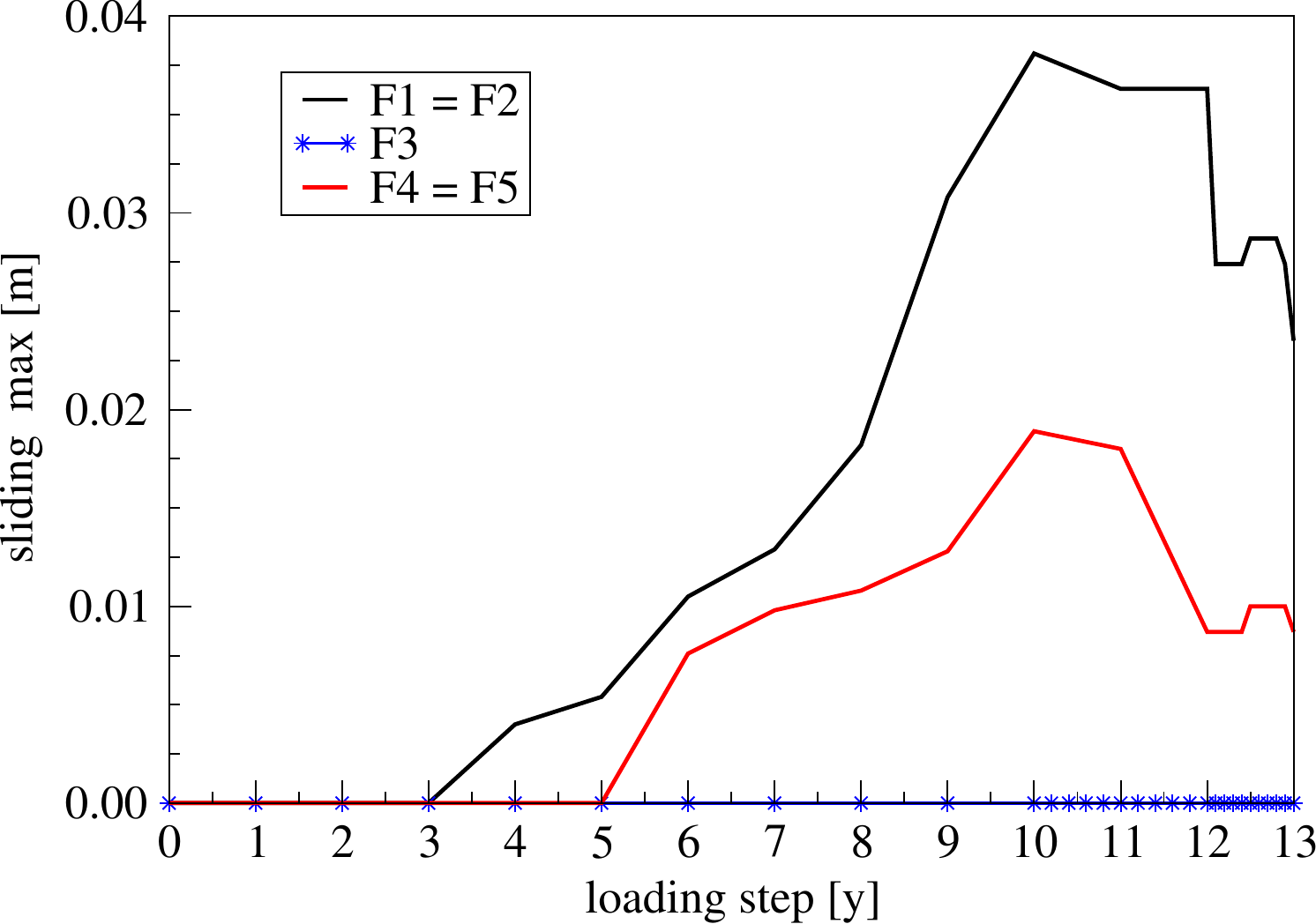}
  \hfill\null
  \caption{Maximum sliding versus time for the investigated scenarios. On the left:
reference case (scenario~1). On the right: using slip-weakening constitutive law
(scenario~2).}
  \label{fig:weak_sliding}
\end{figure}

\subsection{Slip-weakening effect}
The adopted Coulomb frictional criterion can handle slip-weakening effects. Here, the outcome
of a slip-weakening constitutive law for the fault behavior is compared to
that previously obtained using a static friction coefficient equal to $\varphi_s =
30^{\circ}$. The two parameters defining the new constitutive law are $\varphi_d$ and
$D_c$, i.e., the dynamic friction angle and the slip weakening distance, respectively. In
the simulated scenario, the friction angle reduces from $\varphi_s = 30^{\circ}$ to
$\varphi_d = 10^{\circ}$ in a slip distance of $D_c = 2$~mm.

\begin{figure}
  \centering
  \includegraphics[width=\textwidth]{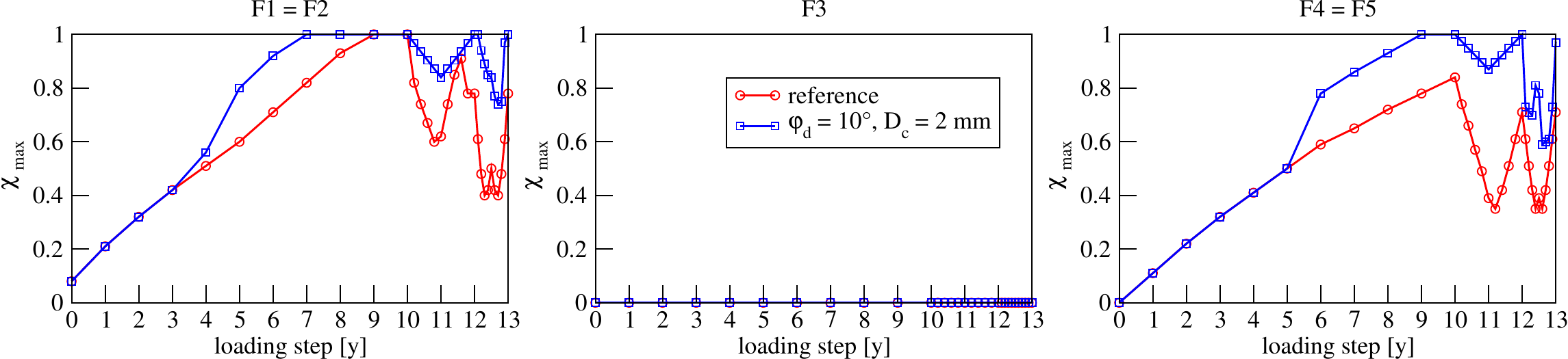}
  \caption{Effect of the Coulomb parameters on $\chi_{\max}$ at increasing loading steps
for each fault. As usual, the pairs F1-F2 and F4-F5 behave identically due to symmetry.
The proposed scenario corresponds to $\varphi_d = 10^{\circ}$ and $D_c = 2$~mm.}
  \label{fig:weakChi}
\end{figure}
\begin{figure}
  \centering
  \includegraphics[width=0.7\textwidth]{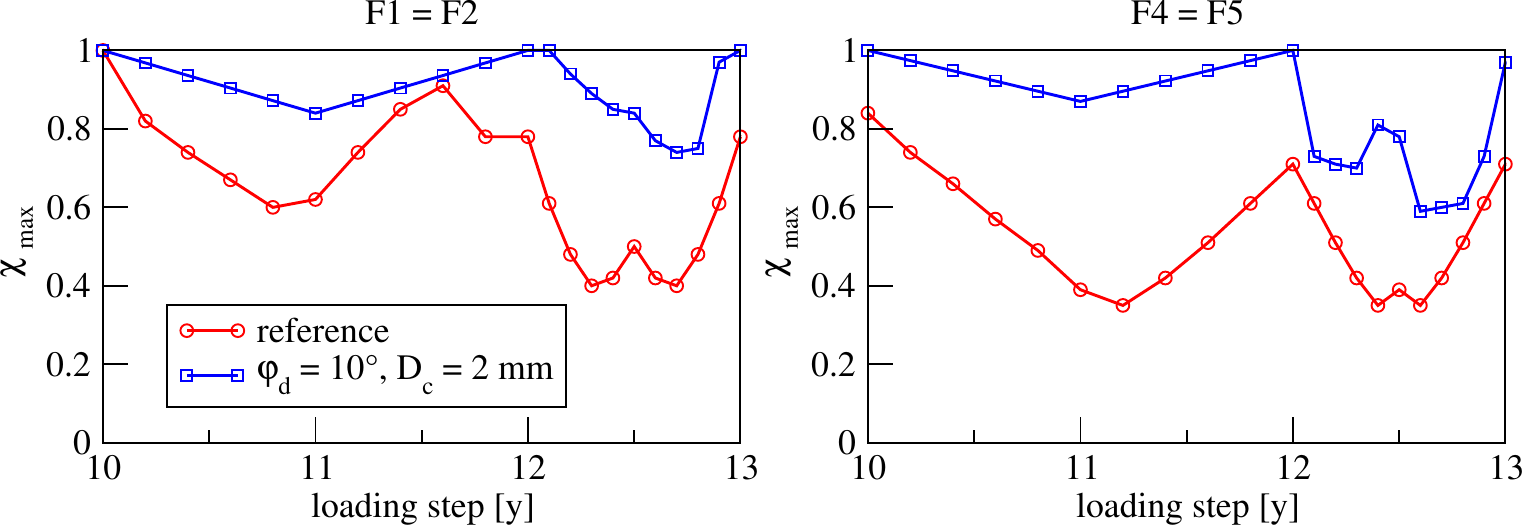}
  \caption{Zoom of Fig.~\ref{fig:weakChi} over the cushion gas injection and UGS phases
for faults F1 = F2 and F4 = F5.}
  \label{fig:weakChiZoom}
\end{figure}

\begin{figure}
  \centering
  \includegraphics[width=0.38\textwidth]{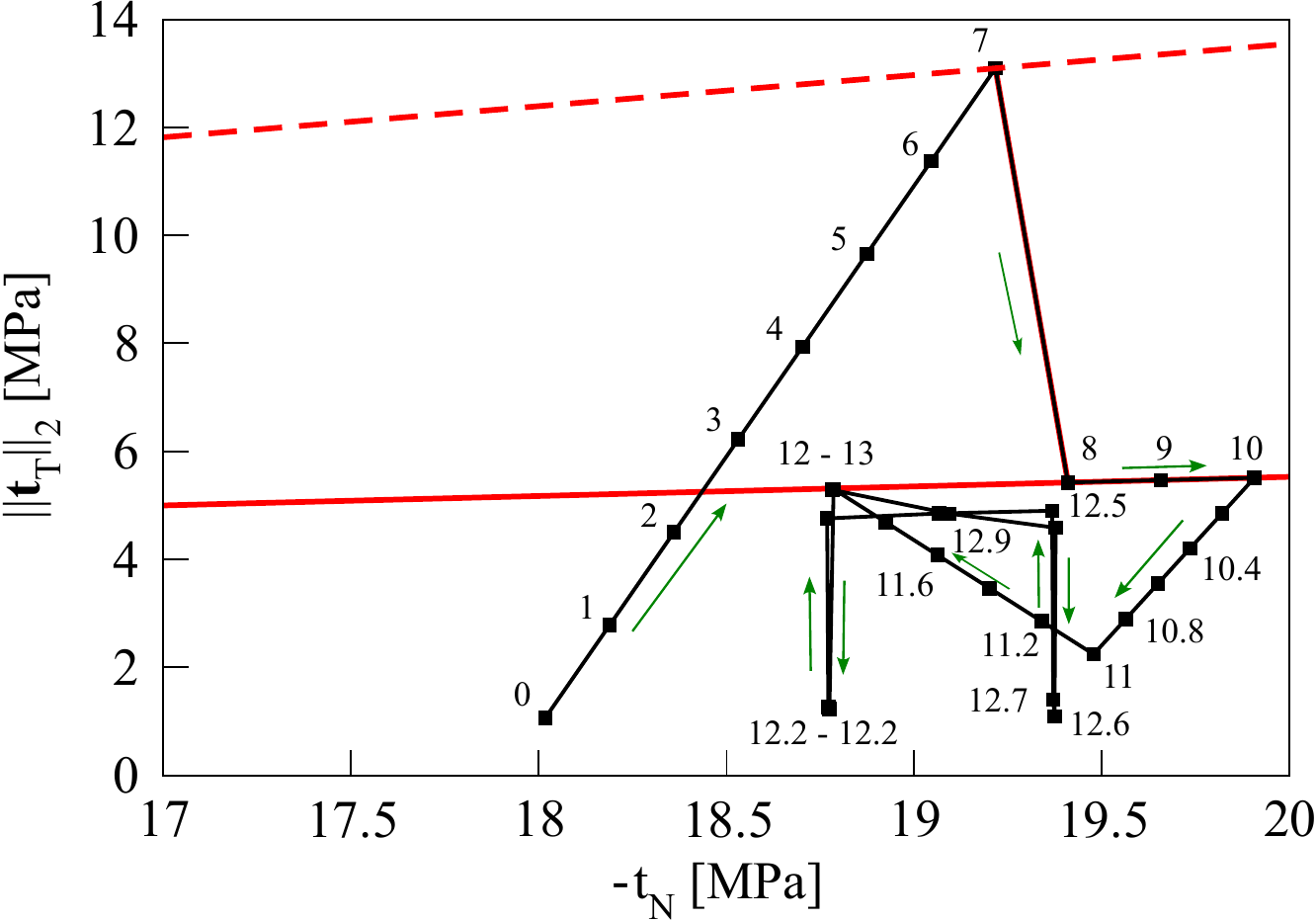}
  \caption{Stress path $\|\vec{t}_T\|_2$ vs $-t_N$ for the F1 element highlighted in
Fig.~\ref{fig:stressPath}. The dashed and the continuous red lines are the yield bound
corresponding to the static condition ($\varphi_s$) and after the slip distance $D_c$ is
overcome, respectively. The numbers along the path denote the loading steps. The primary production (loading steps 1 to 10), the cushion gas injection (loading steps
11 to 12) and the underground gas storage (loading steps 12 to 12.5 -- production -- and 12.5 to 13 -- injection) can be easily recognized.}
  \label{fig:stressPath_w}
\end{figure}

Fig.~\ref{fig:weak_sliding} provides the time behavior of the fault maximum sliding for
the proposed scenario. It can be seen that the current sliding is more than twice that
obtained using a static friction angle. 
Fig.~\ref{fig:weakChi} shows a comparison between the criticality index during the entire
simulation for scenario 1 and 2. It can be noticed that the new constitutive law causes F1, F2, F4 and F5 to slip as well
at the end of the cushion gas and UGS injection phases, but not at loading
step 12.5, i.e., at the end of the 6-month UGS production phase (see the zoom in
Fig.~\ref{fig:weakChiZoom}).
Finally,
Fig.~\ref{fig:stressPath_w} shows the stress path for the same location as in
Fig.~\ref{fig:stressPath}. Because of the reduced friction angle, the yield surface is
reached more easily during PP, at the end of CGI, and at the end of UGS phases. As observed
for the reference scenario, the elastic phases develop with an almost constant normal
stress because of the selected ratios between the reservoir and overburden stiffness and
between the pressure change in the reservoir and within the fault. The stress path and the
yield bound are quite complex due to weakening. Moreover, due to the very small friction
angle ($\varphi = 10^{\circ}$), a large part of UGS is characterized by a the stress state
that develops either on the yield surface or very close to it.

\section{Discussion and conclusions}

\pt{
This two-part study examines the hazards of fault reactivation during underground gas storage (UGS) operations, particularly in the faulted Rotliegend formation in the Netherlands. While seismic events due to pressure increases are well-documented, some ``unexpected'' events occur even when the pressure remains within previously experienced ranges. These events have significant social and economic implications and highlight the need for more reliable predictive tools.}

\pt{In this Part I, a computational model has been implemented to simulate the mechanics of faulted porous rocks. We use Lagrange multipliers to impose the normal and frictional constraints on faults, along with a mixed-dimensional approach and a mixed finite element discretization, where the main unknowns are displacement in the 3D porous body and traction on the fault surfaces. To maintain consistency with classical finite volume discretizations for multiphase flow, we employ low-order hexahedral elements for the 3D continuum and a piecewise constant representation of the traction on the contact surfaces. This approach requires proper stabilization to ensure the regularity of the resulting generalized saddle-point problem. An active-set algorithm and an exact Newton method are implemented to solve the overall nonlinear problem, with ad hoc preconditioning strategies used to facilitate and accelerate the convergence of the inner linear Krylov solver. A discussion of the slip-weakening constitutive law for fault frictional behavior is also included. The model is applied to two realistic scenarios based on a conceptual model derived from an idealization of real UGS fields in the formation of interest. The modeling simulations help identify the key mechanisms that could trigger fault reactivation during UGS activities, even in ``unexpected'' situations where the current stress state seems less demanding than what the porous medium has previously experienced. Using a slip-weakening rheological model for frictional behavior increases the likelihood of fault reactivation during CGI and UGS activities.}

\pt{In Part II, the model will be used to perform an extensive sensitivity analysis to identify key factors influencing fault reactivation. Further work will expand the analysis to include other storage types, such as CO$_2$ and H$_2$ sequestration, and refine the model to define a safe operational bandwidth for UGS management. The ultimate aim is to establish guidelines for managing storage reservoirs and provide a robust framework that can be adapted to other geological settings.}




\section*{Declaration of competing interest}
The authors declare that they have no known competing financial interests or personal
relationships that could have appeared to influence the work reported in this paper.


\section*{Acknowledgements}
This research was supported by the State Supervision of Mines (SodM), Ministry of Economic
Affairs (The Netherlands), Project KEM01 ``Safe Operational Bandwidth of Gas Storage
Reservoirs'' grant.
A.F., M.F., and C.J. from the University of Padova are members of the Gruppo Nazionale
Calcolo Scientifico - Istituto Nazionale di Alta Matematica (GNCS-INdAM).
Computational resources were provided by University of Padova Strategic Research
Infrastructure Grant 2017: ``CAPRI: Calcolo ad Alte Prestazioni per la Ricerca e
l'Innovazione''.

\Urlmuskip=0mu plus 1mu
\bibliographystyle{elsarticle-harv}
\bibliography{biblio}

\clearpage

\end{document}